\documentclass[12pt]{amsart}

\usepackage{amssymb}

\usepackage{amscd}
\usepackage[arrow,matrix,curve,ps]{xy}
\xyoption{dvips}

\newfont{\sheaf}{eusm10 scaled\magstep1}
\numberwithin{equation}{section}

\newcommand{\C}{\ensuremath{\mathbb{C}}}
\newcommand{\CC}{\ensuremath{\mathbb{C}}}

\newcommand{\Z}{\ensuremath{\mathbb{Z}}}

\newcommand{\hol}{\ensuremath{\mathcal{O}}}

\newcommand{\PP}{\ensuremath{\mathbb{P}}}
\newcommand{\Proof}{{\it Proof. }}
\newcommand{\QED}{{\hfill Q.E.D. \medskip}}
\newcommand{\tensor}{{\otimes}}
\newcommand{\sG}{{\mathcal G}}
\newcommand{\sF}{{\mathcal F}}
\newcommand{\sE}{{\mathcal E}}
\newcommand{\sR}{{\mathcal R}}
\newcommand{\sO}{{\mathcal O}}
\newcommand{\sA}{{\mathcal A}}
\newcommand{\sH}{{\mathcal H}}
\newcommand{\sT}{{\mathcal T}}
\newcommand{\sV}{{\mathcal V}}
\newcommand{\Hom}{{\rm Hom}}
\newcommand{\Ext}{{\rm Ext}}
\newcommand{\ra}{{\rightarrow}}
\DeclareMathOperator{\End}{End}

\DeclareMathOperator{\HH}{H}

\DeclareMathOperator{\im}{im}
\DeclareMathOperator{\depth}{depth}

\DeclareMathOperator{\coker}{coker}

\DeclareMathOperator{\Pic}{Pic}

\DeclareMathOperator{\rank}{rank}

\newtheorem{teo}{Theorem}[section]

\newtheorem{lem}[teo]{Lemma}
\newtheorem{cor}[teo]{Corollary}
\newtheorem{rmk}[teo]{Remark}
\newtheorem{prop}[teo]{Proposition}
\newtheorem{slem}[teo]{Sublemma}

\newtheorem{ass}[teo]{Assumption}

\title{ Canonical Projections of Irregular Algebraic Surfaces  }
\author{Fabrizio Catanese and Frank-Olaf Schreyer}

\thanks{
The present research started in 1993 while both authors were visiting 
the Max Planck
Institut f\"ur Mathematik in Bonn. The full results presented here were later
obtained in 1994, when also the second author visited the University 
of Pisa. The
research  continued in the framework of the Schwerpunkt "Globale Methode in der
komplexen Geometrie", and of EAGER. }
\begin{document}
\thispagestyle{empty}
\noindent
{\hfill \today}

\maketitle

We dedicate this article to the memory of Paolo Francia.
\section{Introduction}

  In Chapter VIII of his book
"Le superficie algebriche" F. Enriques raised the question to 
describe concretely the
canonical surfaces with $p_g=4$ and discussed possible constructions of
the regular ones of low degree.

A satisfactory answer to the existence question for degree $\leq 10$ was given
by Ciliberto [1981], and for the case of regular surfaces a 
satisfactory structure
theorem  for the equations of the image surface and its singular locus
was achieved by the first author [1984b].

The main purpose of this paper is to extend those results to the case 
of irregular
surfaces.

Irregular varieties are often easy to construct via transcendental methods,
as is the case for elliptic curves or Abelian varieties. But the 
problem of explicitly
describing the equations of their projective models has always been
a challenge for  algebraic geometry ( cf. Enriques [1949] , Mumford [1966-67]).

 From an algebraic point of view, we might say that irregularity is 
responsible for the
failure of the Cohen-Macaulay condition for the canonical ring of
a variety of general type.

In this context therefore the method of Hilbert resolutions must be replaced
by another method, and we show in this paper that Beilinson's theorem [1978] allows
us to give a suitable generalization of the structure theorem of 
the first author [1984b].

The first two sections are devoted to this extension, and the situation that we consider is the
following: $\varphi$ is a morphism to $\PP^3$ given by four independent sections of the
canonical bundle $K_S$ of a minimal surface of general type, and we assume that the degree of
$\varphi$ is at most two. The following is the main theorem, concerning the case where $\deg \
\varphi =1$.

{\bf Theorem \ref{first main}.}\,{\it
The datum $\varphi : S \rightarrow Y \subset \mathbb P^3 = \mathbb P
(V)$ of a good birational
canonical projection determines a homomorphism
$$(\mathcal O_{\mathbb P} \oplus E)^* (-5) {\buildrel \alpha
\over\longrightarrow}
   (\mathcal O_{\mathbb P} \oplus E),$$
where $E=(K^2 - q+p_g -9) \mathcal O_{\mathbb P}(-2) \oplus q
\Omega_{\mathbb P}^1(-1) \oplus
(p_g - 4)\Omega_{\mathbb P}^2$,
such that
\begin{itemize}
\item[(i)]
$\alpha$ is symmetric,
\item[(ii)]
det $\alpha$ is an irreducible polynomial (defining $Y$),
\item[(iii)]
$\alpha$ satisfies the ring condition (cf. 2.7), and,
\item 
defining $\mathcal F$ as 
the sheaf of $\mathcal O_Y$-algebras given by
the module coker ($\alpha$)
provided with the ring structure determined by  $\alpha$ as in $2.7 (1)$,
\item[(iv)]
{\rm Spec }$\mathcal F$ is a surface with at most rational double points as
singularities.
\end{itemize}

Conversely, given $\alpha$ satisfying i), ii), iii) and iv), $X =
\mbox{\rm Spec}\, \mathcal F$ is the canonical
model of a minimal surface $S$ and
$$\varphi:S \rightarrow X \rightarrow Y \subset \mathbb P^3 = \mathbb P (V)$$
is a good birational canonical projection.}

We recall (cf. \ref{RC}) that the ring condition, or rank condition, on a matrix $\alpha$ is the
condition that the ideal sheaf generated by the top minors of the matrix $\alpha'$ obtained by
deleting the first row of $\alpha$ equals the ideal sheaf generated by the minors of $\alpha$ of
the same size.

The following sections are more in the spirit of Enriques, and we discuss
explicitly the irregular canonical surfaces with $p_g=4$ of lowest degree,
for instance $K^2 = 12$ in the case  of irregularity $q=1$.

We classify completely the above surfaces with   $ p_g = 4, q=1$: they can
be described as a
genus three non hyperelliptic fibration over an elliptic base curve.
 
This  classification shows that the corresponding moduli space has 
only one irreducible component (cf. Theorem 5.10), and we dwell over the geometry of a dense
open set of the moduli space, corresponding to surfaces which we name "of the main stream".

Examples with $p_g = 4, K^2 = 12, q=3$ are the polarizations of type
$(1,1,2)$ on an Abelian threefold: a remarkable subfamily of surfaces for which
the canonical map becomes a degree two covering of a canonical 
surface with $K^2=6$
is given by the "special" surfaces which are the pull backs, under a degree two isogeny, of  the
theta divisor of a principal polarization.

A further example, with $p_g = 4, q=2$ and $K^2 =18$ is provided by certain Abelian covers 
with Galois group $(\Z/2)^2$ of a
principally polarized Abelian surface.

The results presented in this paper were announced in [Catanese, 1997] and very recently
A. Canonaco [2002] was able to extend the method for canonical 
projections to a
weighted $3$-dimensional projective space.

There are still many questions which this paper leaves open:

\begin{itemize}
\item
A precise description of the structure theorem in the general degree two case,
without the assumption that $Y$ be normal
\item
The extension of the structure theorem to the case of higher 
dimensional varieties (cf. however [Catanese 1985] in the "pluriregular case")
\item
The complete classification of canonical surfaces of low degree (this is still
open also in the regular case, as soon as $K^2 \geq 8$, see [Catanese1984b] for the case $K^2
=7$)
\item
The construction of new examples without transcendental methods but 
via computer
algebra.
\item
Decide when is the corresponding moduli space unirational and, in this case, give an explicit
rationally parametrized family
\end{itemize}

\bigskip

\begin{center} NOTATION \end{center}

\bigskip

\noindent
$S$: = the minimal model of a surface of general type \\
${\mathcal R} = {\mathcal R}(S)$ : = the canonical ring of $S$ \\
$X = \mbox{ Proj} \, ({\mathcal R})$ is the canonical model of $S$ \\
$\pi: S \rightarrow X$ is the canonical morphism \\

\noindent
$K$ is a canonical divisor on $X$ or on $S$ (note: $\pi^* (K_X) = K_S)$ \\
$\varphi : S \rightarrow Y \subset \mathbb P^3$ is a good canonical
projection, i.e., \\
$\varphi$ is given by a base point free 4-dimensional subspace $V^*$ of
$\HH^0({\mathcal O}_S (K))$ \\

\noindent
(here $V^*$ denotes the dual vector space of $V$, and ${\mathbb P}^3
= {\mathbb P} (V)$
is the \\
\noindent
projective space of 1-dimensional vector subspaces of $V$). \\
\noindent
$\varphi$ is said to be almost generic if it is good and yields a morphism which is either
birational to
$Y$ or of degree $2$ onto
$Y$ (in [Catanese 1984b] $\varphi$ was said to be quasi-generic if moreover 
$Y$ is  normal in case ${\rm deg}( \varphi ) =2$).
\\

\noindent
$\varphi : S \rightarrow Y$ factors through $\pi$ and a finite morphism
$\psi : X \rightarrow Y$ ;\\

\noindent
$\{ y_0, y_1, y_2, y_3 \}$ : a basis of $V^* \subset \HH^0 ({\mathcal O}_S (K))
\cong \HH^0 ({\mathcal O}_X (K))$ :\\
\noindent
$V^* \oplus W^* \cong \HH^0 ({\mathcal O}_S (K)) \cong \HH^0 
({\mathcal O}_X (K))$,
a fixed splitting; \\
$\sA =$ is the graded polynomial ring $\mathbb C [y_0, y_1, y_2, y_3] =
\mbox{ Sym} \, (V^*)$; \\
$M^{\sim}$ : for a graded $A$-module $M$, denotes the associated sheaf on $\mathbb P^3$; 
therefore  we shall consider 
$${\mathcal R}^{\sim} = \varphi_* {\mathcal O}_S = \psi_*
{\mathcal O}_X =
(\psi_* \omega_X) (-1).$$ 
\noindent
$\sT_Z$ will denote the tangent sheaf of a quasi-projective scheme $Z$

Given Cartier divisors $D, D'$, 

$D \equiv D'$ means that $D$ is linearly equivalent to $D'$, while

$D \sim D'$ means that $D$ is algebraically equivalent to $D'$.

\bigskip
\newpage

\section{Determinantal presentations of canonical projections}

\noindent
By Beilinson's theorem [1978], for any coherent sheaf ${\mathcal F}$
on a projective space $\mathbb P$ there is a complex
$$ \dots \rightarrow {\mathcal C}^i ({\mathcal F}) = \bigoplus_{q-p=i}
\HH^q ({\mathcal F}(-p))
   \otimes_\mathbb C  \Omega_\mathbb P^{p} (p)  \rightarrow {\mathcal C}^{i+1}({\mathcal F})
   \rightarrow  \dots $$
whose cohomology is concentrated in degree 0 and yields ${\mathcal F}$.
For a new proof we also refer to the paper of 
Eisenbud, Fl\o ystad and Schreyer [2001]. In
particular, there is a spectral sequence with $E_1^{-p,q}$ equal to
$\HH^q ({\mathcal F}(-p))
\otimes_\mathbb C  \Omega_\mathbb P^{p} (p)$, and with $d_1^{-p,q}
: E_1^{-p,q}
\rightarrow E_1^{-p+1,q}$ given by the identity tensor

$$id = (\Sigma_j y_j \otimes y_j^*) \in V^* \otimes_\mathbb C
   V = \HH^0({\mathcal O}_\PP (1)) \otimes_\C
\mbox{ Hom}  (\Omega_\PP^1 (1),{\mathcal O}_\PP)$$
which acts according to the tensor rule
$$(x \otimes e) (s \otimes \omega) =
(x s) \otimes ( \omega \neg e),$$
${\neg}$ denoting
contraction of a contravariant tensor
 with a vector .

\medskip
\noindent
 From now on we let ${\mathcal F}(m)$ be $= \psi_* {\sO}_X(m)$, with
$m=0,2$ or $3$ .

 Later on, we shall simply denote ${\mathcal F}(0)$ by ${\mathcal F}$.

The Beilinson table in the particular case $m=3$ reads out as:
$$
\begin{array}{llll}
0 & 0 & 0 & 0 \\
\HH^2({\mathcal O}_X) & \HH^2({\mathcal O}_X (K)) & 0 & 0 \\
\HH^1({\mathcal O}_X) & \HH^1({\mathcal O}_X (K)) & 0 & 0 \\
\HH^0({\mathcal O}_X) & \HH^0({\mathcal O}_X (K)) & \HH^0 ({\mathcal 
O}_X (2K)) & \HH^0
({\mathcal O}_X (3K)) \\
\end{array}
$$
   since $ \HH^3 (\psi_* {\mathcal O}_X (i)) = 0$ for all $i$,
$$ \HH^2 (\psi_* {\mathcal O}_X(i)) = \HH^2 ({\mathcal O}_X (iK))= \HH^0
({\mathcal O}_X (-(i-1)K))^* = 0
\hbox{ for } i \geq 2,$$
and
$$  \HH^1 (\psi_*{\mathcal O}_X (i)) = \HH^1 ({\mathcal O}_X (iK))
= \HH^1({\mathcal O}_X (-(i-1)K))^* = 0 \hbox{ for } i \not= 0,1.$$

\bigskip

\noindent
  From Beilinson's theorem we obtain with a simpler proof  a stronger
version of
a result by Ciliberto

\begin{teo}[Ciliberto, Thm. 2.4,(iii) 1983]
 If $p_g \geq 4$ and $\mid K \mid$ has
no base points, then $\mathcal R$ is generated in degrees $\leq 3$  as an $\sA$-module.
\end{teo}

\noindent
\Proof
The $d_1$-differential on the top row of the Beilinson spectral sequence for
${\mathcal F}(3) = \psi_* {\mathcal O}_X (3)$ has the form
$$
\HH^2({\mathcal O}_X) \otimes_\mathbb C \Omega_\mathbb P^3 (3)=
(W \oplus V)\tensor_\C \Omega_\mathbb P^3 (3) \to
\HH^2({\mathcal O}_X (K)) \otimes_\mathbb C \Omega_\mathbb P^2 (2)
=\Omega_\mathbb P^2 (2)
$$
where the first summand maps to $0$, and the second maps according to the
twisted Serre dual of the Euler sequence: therefore it is surjective with
kernel $\mathcal K$ isomorphic to a direct sum $(W \otimes_\mathbb C \,
{\mathcal O}_\PP (-1))
\oplus \, {\mathcal O}_\PP (-2) \cong [(p_g -4)]{\mathcal O}_\PP (-1) \oplus \,
{\mathcal O}_\PP (-2)$.
\medskip

\noindent
Hence not only ${\mathcal F}(3)$ is a quotient of $\HH^0 ({\mathcal O}_X (3K))
\otimes_\C  {\mathcal O}_\mathbb P$, but ${\sF}(3)$ has
a locally free
resolution by sheaves which are direct sums of sheaves $\sG$ isomorphic
to either ${\mathcal O}_\PP (-2)$, or $\Omega_{\mathbb P}^p (p)$. Such sheaves $\sG$
have the property that
$ \HH^j(\sG(m)) = 0$ for $j > 0$ and $m \ge 0$.

\noindent
Therefore, if we tensor this resolution by ${\mathcal O}_\mathbb
P(m), \, m \geq 0$,
it remains exact on global sections, in particular we have that
$$
\HH^0 ({\mathcal O}_X (3K)) \otimes_\C
\HH^0 ({\mathcal O}_\PP (m-3))
   \rightarrow \HH^0 ({\mathcal O}_X (m  K))$$
is surjective for $m \geq 3$.
\QED

\noindent
In the case $m=2$ the Beilinson table reads out as
$$
\begin{array}{llll}
0 & 0 & 0 & 0 \\
\HH^2({\mathcal O}_X (-K)) & \HH^2({\mathcal O}_X) & \HH^2({\mathcal 
O}_X (K)) & 0 \\
0 & \HH^1({\mathcal O}_X) & \HH^1({\mathcal O}_X (K)) & 0 \\
0 & \HH^0({\mathcal O}_X) & \HH^0({\mathcal O}_X (K)) & 
\HH^0({\mathcal O}_X(2K))
\end{array}
$$
and the symmetry with respect of the middle point of the second row from
the bottom takes each vector space to its Serre dual and each
$d_1$-differential
to its Serre dual map.

\noindent
\begin{slem}
The $d_1$-differential from
$$\HH^0({\mathcal O}_X) \otimes \Omega_\mathbb P^2(2)
\rightarrow \HH^0({\mathcal O}_X(K)) \otimes \Omega_\mathbb P^1(1)$$
   is an isomorphism onto a subbundle.
\end{slem}

\Proof
$d_1$ factors through the natural map of $\Omega_\mathbb P^2(2) \rightarrow
V^* \otimes \Omega_\mathbb P^1(1)$ and the subbundle inclusion of $V^* \otimes
\Omega_\mathbb P^1(1)$ inside $\HH^0({\mathcal O}_X(K)) \otimes
\Omega_\mathbb P^1(1)$,
hence it suffices to show that the first map is a subbundle inclusion. But
this follows from the Beilinson complex for ${\mathcal O}_\mathbb P
(2)$ which yields
the exact sequence
\begin{eqnarray}\label{o2 sequence}
  0 \rightarrow \Omega_\mathbb P^2(2) \rightarrow V^* \otimes
\Omega_\mathbb P^1(1) \rightarrow { \rm Sym}^2  (V^*) \otimes
{\mathcal O}_\mathbb P \rightarrow {\mathcal O}_\mathbb P(2)
\rightarrow 0
\end{eqnarray}
\QED

\noindent
$\psi$ is either birational onto $Y$ or $2:1$ by assumption.
In the second case one can
treat separately the case where $Y$ is a quadric surface (this leads to $K^2_S = 4, q=0$, cf.
Enriques [1949] pages 270-271, i.e., to a double cover of $\PP^1 \times \PP^1$ branched on a
curve of type $(3,3)$). 

Therefore we assume from now on that

\begin{ass}  $ Y$ is not a quadric.
\end{ass}

\noindent
In algebraic terms, this means that
$\HH^0({\mathcal O}_\mathbb P(2))= \mbox{ Sym}^2 \, (V^*)$ is a
direct summand of
$\HH^0({\mathcal O}_X(2K))$, so we can choose a splitting
\begin{eqnarray} \HH^0 ({\mathcal O}_X(2K)) \cong {\rm Sym}^2  (V^*)
\oplus U^*.
   \end{eqnarray}
Moreover, by (\ref{o2 sequence}), we can replace the Beilinson 
complex by a homotopic one,
which gives a new diagram
   $$
\begin{array}{cccc}
0 & 0 & 0 & 0 \cr
\scriptsize \sO(-3) \oplus (U\tensor \sO(-1)) & \scriptsize  W\tensor
\Omega^2(2) & 0 & 0\cr
0 &\scriptsize \HH^1({\mathcal O}_X)\tensor \Omega^2(2)
& \scriptsize \HH^1({\mathcal O}_X (K))\tensor \Omega^1(1) & 0 \cr
0 & 0 &\scriptsize W^* \tensor \Omega^1(1)
& \scriptsize (U^*\tensor \sO )\oplus\sO(2) \cr
\end{array}
$$
Again, here, there is a symmetry taking vector spaces and linear maps to their
Serre duals.

\noindent
Let us denote by $E$ the vector bundle on $\mathbb P^3$ defined by
\begin{eqnarray} \quad E(2) = (U^* \otimes {\mathcal O}_\mathbb P) \oplus
(\HH^1({\mathcal O}_X
(K)) \otimes \Omega_\mathbb P^1(1)) \oplus (W \otimes \Omega_\mathbb
P^2(2))
\end{eqnarray}
We have therefore concluded that ${\mathcal F} = \psi_* {\mathcal
O}_X$ admits a
locally free resolution of length $1$ of the form
\begin{eqnarray}\label{presentation}
  0 \rightarrow ({\mathcal O}_\mathbb P \oplus E)^* (-5)
\rightarrow
{\mathcal O}_\mathbb P \oplus E \rightarrow {\mathcal F} = \psi_*
{\mathcal O}_X
\rightarrow 0.
\end{eqnarray}

Remark that again here, for each twist $m \geq 2$,
(\ref{presentation}) is exact on global
sections.
Now, the two locally free terms of the resolution are dual to each other up
to a twist, and we shall see that one can indeed achieve that the resolution
itself is given by a symmetric map \\

$\alpha : ({\mathcal O}_\mathbb P \oplus E)^*(-5) \rightarrow
{\mathcal O}_\mathbb P
\oplus E \quad (\mbox{that is, } \, \alpha = \alpha^*(-5))$ . \\

\noindent
For simplicity of notation we denote by $\alpha^t$ the map $\alpha^*(-5)$,
and by $F$ the vector bundle ${\mathcal O}_\mathbb P \oplus E$. \\

\noindent
As a first step we state
\begin{lem}\label{lifts} Let ${\mathcal F'}$ be the cokernel 
of another homomorphism
$\beta :F^* (-5) \rightarrow F$.
Then every homomorphism of ${\mathcal F}$ to ${\mathcal F'}$
has a lift to a homomorphism of complexes. Furthermore, any lift of an
isomorphism is an isomorphism of complexes.
\end{lem}

\noindent
\Proof By the exact sequence
$$
\hbox{\rm Hom}(F,F) \rightarrow \hbox{\rm Hom}(F,{\mathcal F'}) \rightarrow
\Ext^1(F, F^*(-5))$$

\noindent
to get a lifting it suffices to have that
$$\Ext^1 (F,F^*(-5)) = \Ext^1(F(2), F^*(-3))=0$$
This holds true since the summands for $F(2), F^*(-3)$ are either $\Omega^j(j)$'s or 
of rank one (and of degree $0$ or $2$ for $F(2)$): moreover, by Bott's formula
$\HH^1(\Omega^j(m)) = 0$,
unless $j=1$ and $m=0$, and by Lemma $2$ of Beilinson [1978],
$\Ext^p(\Omega^j(j),
\Omega^i(i))=0$ for $p \geq 1$.
So every homomorphism of ${\mathcal F}$ to ${\mathcal F'}$ has a lift
to a homomorphism
$f:F \rightarrow F$. Moreover $f$ restricted to $F^*(-5)$ factors through
$g:F^*(-5) \rightarrow F^*(-5)$. This proves the first statement.
\medskip

\noindent
Since also every homomorphism of ${\mathcal F'}$ to ${\mathcal F}$
has a lift by the
same argument, it suffices to prove the second statement for ${\mathcal F'}= {\mathcal
F}$ and for the identity homomorphism of ${\mathcal F}$.
One lift $(f,g)$ is therefore the identity on the complex and this is an
isomorphism. Every other lift $(f_1, g_1)$ of the same automorphism differs
by a homotopy, i.e. $f_1 = f+\alpha \circ h$ and $g_1 = g + h \circ
\alpha$ for some homorphism $h \in \Hom(F, F^*(-5))$ .

\noindent
Since $\HH^0(\Omega_\mathbb P^i(m)) = 0$ for $i \geq1, m \leq i$, and again by
Beilinson's lemma, it follows that
$$
\begin{array}{c}
   \Hom(F,F^* (-5)) = \Hom(F(2), F^*(-3))= \cr
  [ \ \Hom(W,\HH^1({\mathcal O}_X)) \otimes \End(\Omega_\mathbb
P^2(2)) \ ]  \oplus \cr
 \oplus [ \  \Hom(W,W^*)  \otimes \Hom(\Omega_\mathbb P^2(2),
\Omega_\mathbb P^1(1))) \ ] \oplus \cr
\oplus [ \ \Hom(\HH^1({\mathcal O}_x (K), W^*)  \otimes
\End(\Omega_\mathbb P^1(1)) \ ]
\end{array}
$$
On the other hand, if we look at the summands of $\alpha$ involving
$ \Hom(\Omega_\mathbb P^i(i), \Omega_\mathbb P^j(j))$, with $i,j=1$ or $2$,
the only non
zero one is the term in Hom$(\HH^1({\mathcal O}_X) \otimes 
\Omega_\mathbb P^2(2),
\HH^1({\mathcal O}_X(K)) \otimes \Omega_\mathbb P^1(1))$.
Therefore the compositions $\alpha \circ h$ and  $h \circ \alpha$ are
nilpotent and $f_1, g_1$ are isomorphisms, since $f, g$ are  identity matrices.
\QED

\noindent
To obtain a symmetric resolution we apply $\sH
om_{\sO_\PP}(\quad,\omega_\PP(-1))$
to the sequence (\ref{presentation}) and get the exact sequence
\begin{eqnarray} 0 \rightarrow F^* (-5) \rightarrow F \rightarrow
\sE xt^1_{{\mathcal O}_{\mathbb P}}
(\mathcal F, \omega_{\mathbb P}(-1)) \rightarrow 0
\end{eqnarray}
But $\mathcal F (1) = \psi_* \omega_X$, thus by relative duality
for $\psi$, the last term is isomorphic
to $\psi_* \sO_X = \sF$. Pick an isomorphism
$\varepsilon: \mathcal F \rightarrow \sE xt^1_{{\mathcal
O}_{\mathbb P}} (\mathcal F, \omega_{\mathbb P}
(-1))$.
By Lemma \ref{lifts} there is a lift
$$\xymatrix{
0 \ar[r] & F^*(-5) \ar[d]^g \ar[r]^\alpha & F \ar[d]^f\ar[r]
& \sF \ar[d]^\varepsilon \ar[r] & 0 \\
0 \ar[r] & F^*(-5) \ar[r]^{\alpha^t} & F \ar[r] & \sE xt^1_{{\mathcal
O}_{\mathbb P}}
(\mathcal F, \omega_{\mathbb P}(-1))  \ar[r] & 0 \\
}
$$
Applying once more $\sH om_{{\mathcal O}_{\mathbb P}} (\quad
,\omega_{\mathbb P}(-1))$
to the above diagram we get
$$\xymatrix{
0 \ar[r] & F^*(-5) \ar[d]^{f^t} \ar[r]^\alpha & F \ar[d]^{g^t} \ar[r]
& \sE xt^1_{{\mathcal O}_{\mathbb P}}( \sE xt^1_{{\mathcal O}_{\mathbb P}}
(\mathcal F, \omega_{\mathbb P}(-1)),\omega_{\mathbb P}(-1))
\ar[d]^{\varepsilon^\prime} \ar[r] & 0 \\
0 \ar[r] & F^*(-5) \ar[r]^{\alpha^t} & F \ar[r] & \sE xt^1_{{\mathcal
O}_{\mathbb P}}
(\mathcal F, \omega_{\mathbb P}(-1)) \ar[r] & 0 \\
}
$$
and we obtain a canonical isomorphism
\begin{eqnarray}\label{double dual} \sE xt^1_{{\mathcal O}_{\mathbb P}}
(\sE xt^1_{{\mathcal O}_{\mathbb P}}(\mathcal F,
\omega_{\mathbb P}(-1)), \omega_{\mathbb P}(-1)) \cong \mathcal F
\end{eqnarray}
induced by the identity of $F$, since $(\alpha^t)^t = \alpha$. Under
this identification $\varepsilon' =
\sE xt^1_{{\mathcal O}_{\mathbb P}}(\varepsilon, \omega_{\mathbb
P}(-1))$ is the isomorphism
$\mathcal F \rightarrow \sE xt^1_{{\mathcal O}_{\mathbb
P}}(\mathcal F, \omega_{\mathbb P}(-1))$
induced by $g^t$.
\noindent
  Let's assume now that $\varphi$
is birational: then $\Hom_{{\mathcal O}_{\mathbb P}}(\mathcal F, \mathcal
F)= \mathbb C$ and
we have
$\varepsilon'  = \lambda \varepsilon$ for some $ \lambda \in \mathbb \C^*$.
   Moreover by (\ref{double dual})
$$
\begin{array}{rcl}
\varepsilon & = & \sE xt^1_{{\mathcal O}_{\mathbb P}}(\varepsilon',
\omega_{\mathbb P}(-1)),
\mbox{ since both are induced by f,} \\
& = & \sE xt^1_{{\mathcal O}_{\mathbb P}}(\lambda \varepsilon,
\omega_{\mathbb P}(-1))
   = \lambda \, \sE xt^1_{{\mathcal O}_{\mathbb P}} (\varepsilon,
\omega_{\mathbb P}(-1)) \\
&=& \lambda \,
\varepsilon' = \lambda^2 \varepsilon,\\
\end{array}
$$
i.e. $\lambda =\pm 1$ . (We will see in the end that actually
$\lambda = 1.)$ Thus $f$ and $\lambda g^t$
cover the same isomorphism $\varepsilon$, and so does $(f+\lambda
g^t)/2$ . Furthermore
$(f+\lambda g^t)/2$ is an isomorphism by Lemma \ref{lifts}. We claim that
$\beta = ((f+\lambda g^t)/2)
\circ \alpha$ is the desired symmetric matrix. Indeed $\beta$ and
$\alpha$ have isomorphic cokernels,
hence both resolve $\mathcal F = \psi_* \mathcal O_X$ and
$$\beta^t = \alpha^t \circ ((f+\lambda g^t)/2)^t = (\alpha^t \circ f^t
+ \lambda \alpha^t \circ g)/2)
= (g^t \circ \alpha + \lambda f \circ \alpha)/2 = \lambda \beta$$
is either symmetric or skew-symmetric depending on the value of $\lambda$.
Since $\varphi$ is birational  $\det \beta$ gives the equation of
$Y$. In particular $\det \beta$
is not a square. Hence $\beta$ cannot be skew-symmetric, i.e. $\lambda = 1$.

\begin{prop}\label{sym resolution} If $\varphi \colon S \to Y$ is 
birational then
there is a resolution
$$0 \rightarrow (\mathcal O_{\mathbb P} \oplus E)^*(-5) {\buildrel
\alpha \over\longrightarrow}
\mathcal O_{\mathbb P} \oplus E
\rightarrow \mathcal F = \psi_*  \mathcal O_X \rightarrow 0$$
given by a symmetric map $\alpha$ (that is, $\alpha = \alpha^*(-5))$ .
\end{prop}

\noindent
\Proof
Take for $\alpha$ the matrix $\beta$ as above.
\QED

\noindent
\begin{rmk}\label{entries}
The matrix $\alpha$ is a block matrix with entries as indicated below:
$$
\begin{array}{c|cccccccc}
\hbox {Summands} & \sO(-3) & \oplus & U \tensor \sO(-1) &
\oplus &\HH^{0,1} \tensor \Omega^2(2) & \oplus& W^* \tensor \Omega^1(1) \\
\hline
&&&&&&&\\
\sO(2) & S_5 V^* && S_3 V^* && \HH^0(\Lambda^2 \sT_\PP) && \HH^0(\sT_\PP) \\
\oplus &&&&&&&\\
U^* \tensor \sO & S_3 V^* && V^*\cong \Lambda^3 V && \Lambda^2 V && V \\
\oplus &&&&&&&\\
\HH^{2,1}\tensor \Omega^1(1)& \HH^0(\Lambda^2 \sT_\PP)&&\Lambda^2V&&V &&0\\
\oplus &&&&&&&\\
W\tensor\Omega^2(2)&\HH^0(\sT_\PP)&&V&&0&&0\\
\end{array}
$$
\end{rmk}

\noindent
Note that the $h^{2,1} \times h^{0,1}$ block is actually
skew-symmetric, since it is induced from
wedge-product
$$\HH^1(\mathcal O_X) \times \HH^1(\mathcal O_X) \rightarrow \HH^2(\mathcal
O_X) = \HH^0 (\sO_X(K))^*$$
composed with the projection
$$\HH^0(\mathcal O_X(K))^* \rightarrow \HH^0(\mathcal O_{\mathbb P}(1))^* = V$$
However each element of
$$V \cong \mbox{Hom}(\Omega_{\mathbb P}^2(2), \Omega_{\mathbb P}^1(1)
\cong \mbox{Hom}((\Omega_{\mathbb P}
^1(1))^*, (\Omega_{\mathbb P}^2(2))^*)$$
\noindent
gives a skew-symmetric morphism of bundles. So the resulting morphism\\
$$\HH^{0,1} \otimes \Omega_{\mathbb P}
^2(2) \rightarrow \HH^{2,1} \otimes \Omega_{\mathbb P}^1(1)$$
is symmetric again. For general sign patterns in Beilinson monads
of symmetric or skew-symmetric sheaves see Eisenbud and Schreyer [2001].
\medskip

\noindent
For a morphism $\alpha: G \rightarrow F$ of vector bundles on a scheme
$Z$ we denote by $I_r(\alpha)$
the ideal sheaf of  $r\times r$ minors of $\alpha$, i.e. the image of
$\Lambda^r G \otimes (\Lambda^r F)^*
\rightarrow \mathcal O := \mathcal O_Z$ under the natural map induced by
$\Lambda^r(\alpha)$. So $I_r(\alpha)$ is the
$(\rank(F)-r)^{\mbox{th}}$ Fitting ideal of coker$(\alpha)$.

\begin{teo}\label{RC}[Catanese {[1984b]}, de Jong and van Straten {[1990]}]
Let $\alpha = (\alpha_1, \alpha'): G \rightarrow \mathcal O \oplus E$
be a morphism of vector bundles
with $r= \rank E = \rank G-1$. Suppose $\det(\alpha)$ 
is  a  non zero-divisor and ($(\det(\alpha))= I_{r+1}(\alpha)$)
 $\depth(I_r(\alpha'),\mathcal O)
\geq 2$ . Let  $Y \subset Z$ denote the
subscheme defined by $\det(\alpha)$.
Then the following are equivalent:
\begin{itemize}
\item[(1)] $\mathcal F = \coker(\alpha)$ carries the structure of a sheaf
of commutative $\mathcal O_Y$-algebras
with $1 \in \sF$ given by the image of $1 \in \Gamma(Z,\mathcal O)
\subset \Gamma(Z,\mathcal O \oplus E)$,
\item[(R.C.)] $I_r(\alpha) = I_r(\alpha')$ .
\end{itemize}
\end{teo}

\noindent
\Proof
Since $\det(\alpha)$ is a non-zero divisor,
$$
0 \to  G {\buildrel \alpha \over \longrightarrow} F \to \sF \to0
$$
with $F= \mathcal O \oplus E$ is exact. As an $\mathcal O_Y$-module
$\mathcal F$ has an infinite
periodic resolution
$$\ldots \rightarrow \mathcal B^2 \otimes F_Y  {\buildrel \beta_Y
\over \longrightarrow}
   \mathcal B \otimes
G_Y {\buildrel \alpha_Y \over \longrightarrow} \mathcal B \otimes
F_Y {\buildrel \beta_Y \over \longrightarrow} G_Y {\buildrel \alpha_Y
\over \longrightarrow}  F_Y \rightarrow \mathcal F \rightarrow 0$$
where $\mathcal B = \Lambda^{r+1}G \otimes (\Lambda^{r+1}F)^*$,
$-_Y = -\otimes \mathcal O_Y$ denotes
restriction to $Y$, and the map
$$\beta: \Lambda^{r+1}G \otimes \Lambda^{r+1} F^*\tensor F \rightarrow G$$
is induced by $\Lambda^r(\alpha)$, i.e. $\beta$ is given by the
matrix of cofactors of $\alpha$. Exactness
follows, since $\alpha \cdot \beta = \det (\alpha)\mbox
{id}_F$ and $\beta \cdot \alpha = \det
(\alpha)\mbox {id}_G$ and det$(\alpha)$ is a non-zerodivisor, by
Eisenbud [1980]. 

The above resolution is obtained by truncating the infinite exact periodic complex

$$\ldots \rightarrow \mathcal B^2 \otimes F_Y  {\buildrel \beta_Y
\over \longrightarrow}
   \mathcal B \otimes
G_Y {\buildrel \alpha_Y \over \longrightarrow} \mathcal B \otimes
F_Y {\buildrel \beta_Y \over \longrightarrow} G_Y {\buildrel \alpha_Y
\over \longrightarrow}  F_Y {\buildrel \beta_Y
\over \longrightarrow}\mathcal B^{-1} \otimes
G_Y {\buildrel \alpha_Y \over \longrightarrow} \mathcal B^{-1} \otimes
F_Y  \rightarrow ...$$

Similarly we have an exact
infinite periodic complex
$$\ldots  \longrightarrow \mathcal B \otimes
F^*_Y {\buildrel \alpha_Y^t \over \longrightarrow}
   \mathcal B \otimes
G^*_Y {\buildrel \beta_Y^t \over \longrightarrow}
F^*_Y {\buildrel \alpha_Y^t \over \longrightarrow} G^*_Y {\buildrel
\beta_Y^t \over \longrightarrow} \mathcal B^{-1} F^*_Y
\longrightarrow \ldots$$
\noindent
It follows  then that $\mathcal C:=\sH om_Y(\mathcal
F,\mathcal O_Y)$ satisfies $\mathcal F
= \sH om_Y(\mathcal C, \mathcal O_Y)$.

\noindent
Recall that $F= \mathcal O \oplus E$ whence $\mathcal O_Y \subset
\mathcal F$ by Cramer's rule. By the assumption
$\depth(I_r(\alpha'), \mathcal O) \geq 2$ the quotient $\mathcal
F/\mathcal O_Y$, which is annihilated by
$I_r(\alpha')$, satisfies $\sH om_Y(\mathcal F/\mathcal O_Y,
\mathcal O_Y)=0$.
Therefore
$$
\begin{array}{lcl}
\mathcal O_Y & = & \sH om_Y (\mathcal O_Y, \mathcal O_Y) \supset
\mathcal C = \sH om_Y(\mathcal F,
\mathcal O_Y)\\
& \cong & \ker (F^*_Y\rightarrow G^*_Y) \cong \im (\mathcal
B \otimes G^*_Y \rightarrow F^*_Y), \\
\end{array}
$$
so $\mathcal C \subset \mathcal O_Y$ is the ideal sheaf
$$
\mbox{im}(\mathcal B \otimes G^*_Y \rightarrow \mathcal O_Y) =
I_r(\alpha')/(\mbox{det}(\alpha)).
$$
\noindent
Suppose $\mathcal F$ is a ring. Then $\mathcal C \subset \mathcal O_Y$
is called the conductor of
$\mathcal O_Y \subset \mathcal F$ and $\mathcal F \subset
\sH om_Y(\mathcal C, \mathcal C)$, i.e.
the image of $\mathcal C \times \mathcal F$ in $\mathcal O_Y$ is
contained in $\mathcal C \subset
\mathcal O_Y$: in fact for each $(c,m)$ in $\mathcal C \times
\mathcal F$ the image $c(m)\in \sO_Y$ carries $\mathcal F$
into $\mathcal O_Y$ since $c(m)\mathcal F = c(m \mathcal F) \subset
\mathcal O_Y$. Thus
$$\mathcal F \subset \sH om_Y(\mathcal C,\mathcal C) \subset
\sH om_Y(\mathcal C, \mathcal O_Y)
= \mathcal F.$$
In particular $\sH om_Y(\mathcal C,\mathcal C) =
\sH om_Y(\mathcal C,\mathcal O_Y).$
Now the last equality means that every $\varphi \in
\sH om_Y(\mathcal C, \mathcal O_Y) = \mathcal F =
\coker(G_Y \rightarrow F_Y)$ has image in $\mathcal C =
I_r(\alpha')/(\mbox{det}(\alpha))$.
That is the pairing induced by $\beta_Y$
$$F_Y \times (\mathcal B \otimes G^*_Y) \rightarrow \mathcal O_Y$$
has image in $I_r(\alpha')\mathcal O_Y$. Since $\beta$ is the matrix
of cofactors of $\alpha$,
$\sH om_Y(\mathcal C, \mathcal C) = \sH om_Y(\mathcal C,
\mathcal O_Y)$ is equivalent to
$I_r(\alpha) = I_r(\alpha')$. \medskip

\noindent
Conversely, if $(R.C.)$ holds, we have $\mathcal F =
\sH om_Y(\mathcal C, \mathcal O_Y)  =
\sH om_Y(\mathcal C, \mathcal C)$ is an $\mathcal O_Y$-algebra under
composition. The structure is commutative, since  $\mathcal C$ as
$\mathcal O_Y$-module is invertible
on $Y - V(I_r(\alpha'))$ and $\depth(I_r(\alpha'), \mathcal F) \geq 1$
by assumption.

\QED

\begin{rmk}\label{further RC}
We call the condition
${\rm (R.C.)}\,  I_r(\alpha)=I_r(\alpha')$ the Ring
Condition (or Rank Condition). For a
symmetric matrix
$$\alpha=\begin{pmatrix}
\alpha_{11} & \alpha_{12} \cr
\alpha_{12}^t &\alpha'' \cr
\end{pmatrix}
:(\mathcal O \oplus E)^* (-5) \rightarrow \mathcal O \oplus E$$
the ring condition implies the further rank condition
$$I_{r-1}(\alpha') = I_{r-1}(\alpha'')$$
for lower minors of $(\alpha')^t = (\alpha_{12},\alpha''): E^*(-5)
\rightarrow \mathcal O \oplus E$
in case depth$(I_{r-1}(\alpha''),\mathcal O) = 3$ has the expected
maximal value.
\end{rmk}
Cf. Prop. 4.1 of
Mond and Pellikaan [1987], and  Prop. 5.8 and 5.10 of Catanese [1984b].

\begin{teo}\label{first main}
The datum $\varphi : S \rightarrow Y \subset \mathbb P^3 = \mathbb P
(V)$ of a good birational
canonical projection determines a morphism
$$(\mathcal O_{\mathbb P} \oplus E)^* (-5) {\buildrel \alpha
\over\longrightarrow}
   (\mathcal O_{\mathbb P} \oplus E),$$
where $E=(K^2 - q+p_g -9) \mathcal O_{\mathbb P}(-2) \oplus q
\Omega_{\mathbb P}^1(-1) \oplus
(p_g - 4)\Omega_{\mathbb P}^2$,
such that
\begin{itemize}
\item[(i)]
$\alpha$ is symmetric,
\item[(ii)]
det $\alpha$ is an irreducible polynomial (defining $Y$),
\item[(iii)]
$\alpha$ satisfies the ring condition, and
\item
defining $\mathcal F$ as the sheaf of $\mathcal O_Y$-algebras given by
the module coker ($\alpha$)
provided with the ring structure determined by  $\alpha$ as in $2.7(1)$,
\item[(iv)]
{\rm Spec } $\mathcal F$ is a surface with at most rational double points as
singularities.
\end{itemize}

Conversely, given $\alpha$ satisfying i), ii), iii) and iv),  $X =
\mbox{\rm Spec}\, \mathcal F$ is the canonical
model of a minimal surface $S$ and
$$\varphi:S \rightarrow X \rightarrow Y \subset \mathbb P^3 = \mathbb P (V)$$
is a good birational canonical projection.
\end{teo}

\noindent
\proof
The first statement follows by combining \ref{RC} and \ref{sym resolution}.
Notice
that, since $\det\alpha$ is irreducible (in one direction, this is a consequence of the
birationality of $\varphi$, in the other direction, it holds by assumption), the ring condition
$I_r(\alpha)=I_r(\alpha')$ implies depth$(I_r(\alpha'),\mathcal O_{\mathbb P}) \geq 2$, since
otherwise all $r\times r$ minors of $\alpha$ would have a common
irreducible factor, whose square would
divide $\det(\alpha)$. So the assumptions of \ref{RC} are satisfied.

\noindent
For the converse, we note that duality
applied to
$$\psi:X=\mbox{Spec}\, \mathcal F \rightarrow Y \subset \mathbb P^3$$
gives
$$\mathcal F(1) = \psi_* \omega_X,$$
since $\alpha$ is symmetric. So $\psi^* \mathcal O_{\mathbb P} (1) \cong
\omega_X$, and, because $X$ has only rational
double points as singularities, $\varphi^* \mathcal O_{\mathbb P} (1) \cong
\omega_S$ holds on the desingularization
$S$ of $X$. So $\varphi:S \rightarrow Y \subset \mathbb P^3$ is a
quasi generic birational canonical
projection.

\QED

\begin{rmk}
Given a symmetric matrix
$\alpha:(\mathcal O \oplus E)^* (-5) \rightarrow \mathcal O \oplus E$
as in the theorem, we denote by
$\alpha' : (\mathcal O \oplus E)^* (-5) \ra  E$ and $\alpha'' :E^* (-5)
\rightarrow E$
the distinguished submatrices. Then $\det(\alpha'')$ defines the
adjoint surface $F$ of degree $K^2-5$.
$F$ intersects $Y$ precisely in the non-normal locus $\Gamma$ of $Y$,
which is defined by $I_r(\alpha')$.
$F$ is singular at the points of the subscheme $T$ defined by $I_{r-1}(\alpha'')$,
$\Gamma$ has embedding dimension
$3$ at $T$, and the points of $T$ are at least triple for $Y$. \end{rmk}

  Typically the points of $T$ correspond to triple points of
$Y$, and $F$ has ordinary quadratic
singularities in $T$.
In terms of the invariants $d=K^2, q$ and $p_g$ of $S$ we have
$$\deg \Gamma = 1/2 \ d^2-5/2 \ d+1$$
and
$$\deg T = 1/6 \ d^3-5/2 \ d^2+37/3 \ d-4(1-q+p_g).$$

\begin{prop}\label{albanense curve bound} Suppose 
$\varphi : S \to Y \subset \PP^3=\PP(V)$ 
is a good canonical projection with $Y$ not a quadric. 
If the Albanese image of $S$ is a curve
then
$$ K^2 \ge 4q+4+2(p_g-4) = 2p_g + 4q -4.$$
If equality holds then the map is not birational.
\end{prop}

\proof
Since $ Y \subset \PP^3$ is a surface we have the  presentation \ref{presentation}
of $\phi_* \sO_S=\psi_* \sO_X$ given by a matrix $\alpha$ with entries
as indicated in Remark \ref{entries} (with $\alpha$ perhaps not 
symmetric). If the Albanese dimension is one then
$$\HH^0(\Omega_S^1) \times \HH^0(\Omega^1_S) \to \HH^0(\Omega^2_S)$$
is the zero map, and so is $\HH^1(\sO_S)\times\HH^1(\sO_S) \to \HH^2(\sO_S)$
by Hodge symmetry. Thus we have a large block of zeroes. On the other
hand the determinant of $\alpha$ equals the equation of $Y$ to the power
$\deg(\varphi)$, in particular $\det \alpha \not= 0$. Thus the $(q+p_g-4) \times
(q+p_g-4)$ block  cannot be too big. More precisely,
$$ 1+ \dim U \ge 3(q+p_g-4).$$
Since $1+\dim U = 1+ K^2 + 1-q+p_g - \dim S_2 V = K^2  -q + p_g -8$  the desired
inequality follows. Moreover in case of equality we have
$$\det \alpha =\lambda [\det( \sO(-3) \oplus (U \tensor \sO(-1)) 
\to (\HH^{2,1} \tensor\, \Omega^1(1) )\oplus (W \tensor\, \Omega^2(2))]^2$$
for some scalar $\lambda \in \CC$. Thus $\varphi $ cannot be birational.
\QED

\begin{rmk}
The same argument, but without the assumption that the Albanese image of
$S$ be a curve, gives in general the inequality $ K^2 \geq 2 p_g - 2q -4$ which is
however weaker than Noether's inequality $ K^2 \geq 2 p_g  -4$. For $ q \geq 1$ we may get 
$ K^2 \geq 2 p_g - 2q + 2$, which is still however weaker than the
inequality given by Debarre [1982] for irregular surfaces, namely, $ K^2 \geq 2 p_g $.
\end{rmk}

\section{The case of double covers}

\noindent
Suppose the good canonical projection $\psi:S
\rightarrow Y \subset \mathbb P^3$
is $2:1$ and that $Y$ is not a quadric. Then $\mathcal F = \psi_*
\mathcal O_X$ still has a locally
free resolution of length $1$ of the form (\ref{presentation})
$$ \quad 0 \rightarrow (\mathcal O_{\mathbb P} \oplus E)^*(-5)
\rightarrow \mathcal O_{\PP}
\oplus E \rightarrow \mathcal F = \psi_* {\mathcal O}_X \rightarrow 0.$$
with $E$ the vector bundle on $\mathbb P^3$ defined by
$$ E(2) = (U^* \otimes \mathcal O_{\mathbb P}) \oplus
(\HH^1(\mathcal O_X(K)) \otimes \Omega_{\mathbb P}^1(1))
\oplus (W \otimes \Omega_{\mathbb P}^2(2)).$$
However the proof that the resolution can be chosen symmetric needs
further arguments.
Let $\rho :Z \rightarrow Y$ denote the normalization of $Y$. Then
$\psi : X \rightarrow Y$ factors
over $\varepsilon : X \rightarrow Z, \hbox{ where } \varepsilon$ is $2 : 1$.
The
covering involution $\Phi:
X \rightarrow X$ induces a decomposition of $\varepsilon_*
\mathcal O_X$ into invariant and
anti-invariants parts:
\begin{eqnarray}
\varepsilon_* \mathcal O_X = \mathcal O_Z \oplus
\mathcal H.
\end{eqnarray}
\noindent
On $Y$ this induces the decomposition
\begin{eqnarray}
   \mathcal F = \psi_* \mathcal O_X = \rho_* \mathcal
O_Z \oplus \rho_* \mathcal H.
\end{eqnarray}
This in turn decomposes the Beilinson cohomology groups of $\mathcal
F$ and (assuming that $Y$ is not
a quadric), this gives a decomposition of (\ref{presentation}). There are two
cases how the isomorphism
$\mathcal F \rightarrow {\sE xt}^1_{\mathcal O_{\mathbb P}}
(\mathcal F, \omega_{\mathbb P}(-1))$
can respect the summands (which are generically of rank $1$ on $Y$).
Either
\begin{itemize}
\item[(a)]     $\quad \rho_* \mathcal O_Z \cong {\sE xt}^1_{\mathcal
O_{\mathbb P}}
(\rho_* \sO_Z, \omega_{\mathbb P}(-1))$
\end{itemize}
or
\begin{itemize}
\item[(b)]     $\quad \rho_* \mathcal O_Z \cong
{\sE xt}^1_{\mathcal O_{\mathbb P}}
(\rho_* \sH, \omega_{\mathbb P}(-1))$
\end{itemize}
Case (a) occurs when  $y_0, y_1, y_2, y_3$
pullback to $\Phi$-invariant sections of $\HH^0(\sO_X(K))$
i.e. $V^* \subset \HH^0(\sO_X(K))^+$.

 Case (b) occurs when 
 $V^* \subset \HH^0(\sO_X(K))^-$.

\noindent
{\bf The $\Phi$-invariant case (a).}

Since   $V^* \subset \HH^0(\sO_X(K))^+$, the isomorphism 
$\psi_* \omega_X$ with $\psi_* \hol_X(1)$ respects the
invariant and anti-invariant summands.
Therefore, $$ (\psi_* \omega_X)^+ = \rho_* \omega_Z = {\sE xt}^1_{\mathcal
O_{\mathbb P}}
(\rho_* \sO_Z, \omega_{\mathbb P})$$ is isomorphic to 
$\rho_* \hol_Z (1)$ as asserted. Since moreover $\rho$ is birational,
$\rho^* \rho_* \hol_Z (1) = \hol_Z (1) $ and by the projection formula we infer
that $\sO_Z(1)\cong \omega_Z$, in particular, that $Z$ is Gorenstein, whence
the canonical model of a surface of general type.

Similarly, we see that $ (\psi_* \omega_X)^- = {\sE xt}^1_{\mathcal O_{\mathbb P}}
(\rho_* \sH, \omega_{\mathbb P}) = \rho_* \sH (1).$

Decomposing then
\begin{eqnarray} E = E_+\oplus E_- \end{eqnarray}
into parts coming from $ \rho_* \mathcal
O_Z \oplus \rho_* \mathcal H$ we obtain in  this case  a
decompositon of (\ref{presentation}) into

\noindent
\begin{eqnarray}\label{sequence +}
   0 \rightarrow  (\sO_\PP \oplus E_+)^* (-5)
{\buildrel \alpha_+ \over \longrightarrow}
    \sO_\PP \oplus E_+ \rightarrow \rho_*
  \sO_Z
\rightarrow 0 \end{eqnarray}
and

\noindent
\begin{eqnarray}\label{sequence -} 0 \rightarrow  (E_-)^* (-5)
{\buildrel \alpha_- \over \longrightarrow}
   E_- \rightarrow \rho_*
\sH
\rightarrow 0 \end{eqnarray}

\noindent
The argument of \ref{sym resolution} shows
that both $\alpha_+$ and $\alpha_-$  can be
   chosen to be symmetric matrices. 
 We can also apply Theorem \ref{first main} to $Z$,
since  as we observed $Z = X/\Phi$ is a
 canonical model of a surface.  (\ref{sequence +} is the determinantal description
ot the good birational  canonical projection $\rho: Z \to Y \subset
\PP^3=\PP(V).$

We leave aside for the time being the task of describing the  $\sO_Z$ module
structure on $\rho_*\sH$: we simply observe that the bilinear map 
$$(E_-) \times (E_-) \ra \hol_{\PP}$$ is given as in [Catanese 1981] by the adjoint
matrix of $\alpha_-$, which solves the problem in the very particular
case where $Y$ is normal.

\noindent
{\bf  The $\Phi$-anti-invariant case (b). }
Here (\ref{presentation}) decomposes into

\noindent
\begin{eqnarray}\label{sequence +-} 0 \rightarrow ( E_-)^* (-5)
{\buildrel \alpha_+ \over \longrightarrow}
{\mathcal O}_\mathbb P \oplus E_+ \rightarrow \rho_*
{\mathcal O}_Z
\rightarrow 0 \end{eqnarray}
and

\noindent
\begin{eqnarray} 0 \rightarrow ({\mathcal O}_\mathbb P \oplus E_+)^* (-5)
{\buildrel \alpha_- \over \longrightarrow}
   E_- \rightarrow \rho_* \sH
\rightarrow 0 \end{eqnarray}

\noindent
and we may choose $\alpha_-=(\alpha_+)^t$ . Notice that $\alpha_+$
satifies the ring condition (R.C.).

\medskip \noindent
To recover $X$ from (\ref{sequence +}) and (\ref{sequence -}) or
from (\ref{sequence +-}) we need in addition to describe the map

\noindent
\begin{eqnarray}   S_2(\rho_* \sH) \to \rho_* \sO_Z. \end{eqnarray}

\section{ Generalities on irregular surfaces with $p_g= 4$ .}

As explained in the introduction, one of the main purposes of our investigation
is to understand the equations of the projections of irregular canonical surfaces
in $\PP^3$.

For this reason the most natural case to consider is the case where $p_g =4$,
and  there is no choice whatsoever to make for the projection. 

We recall, for the reader's benefit, some important inequalities for irregular surfaces
\begin{itemize}
\item
Castelnuovo' s Theorem [1905] : $\chi(S) \geq 1$ if the surface $S$ is not
ruled.
\end{itemize}
From Castelnuovo's theorem follows that if $p_g=4$, then  the irregularity $q(S)$ is
$\leq 4$.

We have also the 
\begin{itemize}
\item
Inequality of Castelnuovo-Beauville (Beauville [1982]):
$ p_g \geq 2q -4 $, equality holding if and only if $S$ is a product of
a curve of genus $2$ with a curve of genus $g \geq 2$.
\end{itemize}

Whence follows that, if $p_g = q =4$, then $S$ is the product of two curves of
genus $2$ and its canonical map is a $(\Z/2)^2$-Galois covering of a smooth
quadric.

Debarre has moreover shown [1982] that for an irregular surface one has the
following 
\begin{itemize}
\item
Debarre's inequality: if $S$ is irregular, then $K^2 \geq 2 p_g$.
\end{itemize}
Therefore, in our case, the above inequality yields more than the more general
inequalities by Castelnuovo [1891] and  by Horikawa [1976b]], Reid [1979], Beauville
[1979] and Debarre  [1982]:   $ K^2 \geq 3 p_g + q -7 $ if the
canonical map is birational.

Our inequality $ K^2 \geq 2 p_g + 4 q - 4 $ if the Albanese image is a curve
severely restricts the numerical possibilities if $p_g=4 , q=3$, since by
the Bogomolov-Miyaoka-Yau inequality we always have $ K^2 \leq 9 \chi$, thus
$ 16 \leq K^2 \leq 18$ if the Albanese map is a pencil. 

This case is completely solved by using the following inequalities for surfaces fibred over
curves:

\begin{itemize}
\item
Arakelov's inequality: let $f : S \ra B $ be a fibration onto a curve $B$ of genus $b$,
with fibres curves of genus $g \geq 2$. Then $K^2_S \geq 8 (b-1) (g-1)$, equality holding 
only if the fibration has constant moduli.
\end{itemize}

\begin{itemize}
\item
Beauville's inequality: let $f : S \ra B $ be a fibration onto a curve $B$ of genus $b$,
with fibres curves of genus $g \geq 2$. Then $\chi(S) \geq  (b-1) (g-1)$, equality holding 
if and only if the fibration is an etale bundle (there is an etale cover $B' \ra B$ such that the
pull back is a product $B' \times F$).
\end{itemize}

By Beauville's inequality follows that if $q=3$ and the Albanese image is a curve, then (take
as $B$ the genus $3$ curve which is the Albanese image,and $f$ the Albanese map)  the
Albanese map is an etale bundle with fibre $F$ of genus $2$.

In particular, we have $K^2_S = 16$ and all our surfaces are obtained as follows:

let $G$ be a finite group acting faithfully on a curve $F$ of genus $2$ in such a way that $F /G
\cong \PP^1$, and take an etale $G$-cover $B' \ra B$ of a curve $B$ of genus $3$:
then our surfaces $S$ with $p_g= 4$, $q=3$, $K^2_S = 16$ are exactly the quotients 
$(F \times B') / G$ of the product $(F \times B')$ by the diagonal action of $G$.

The groups $G$ as above were classified by Bolza [1888] (cf. also [Zucconi1994]). 

If instead the Albanese image is a curve of genus $q=2$, then, since we assume $p_g=4$, then
$\chi(S) = 3 \geq (g-1)$, and the genus of the Albanese fibres is at most $4$.

The case $g=4$ gives again rise to an etale bundle with fibre $F$ of genus $4$, so that our
surface $S$ is a quotient $ S = F \times B' / G$ where $ F/ G \cong \PP^1$ and $  B =
B'/G$ is a curve of genus $b=2$ (one must impose the condition that $G$ operate freely on the
product).

In particular, we have $K^2_S = 24$. A concrete example is furnished by $G = \Z / 2$, which
operates freely on the genus $3$ curve $B'$. Then  $ H^0(K_S) = H^0(K_{F\times B'})^+
= H^0(K_F)^- \otimes H^0(K_{B'})^- $, and since $ H^0(K_{B'})^-$ is $1$-dimensional, the
canonical system is a pencil and the canonical image of $S$ is the canonical image of $F$,
namely,  a twisted cubic curve in
$\PP^3$.

 It would be
also interesting to determine what is the minimal value of
$K^2$ for an irregular surface with $p_g=4$ such that the canonical map is birational.

In the next section we shall see that the answer to this question is: 
$K^2 = 12$ in the case $q=1$, and in this case we shall give a complete classification 
of the surfaces for which the canonical map is birational.

On the other hand, the problem remains for $q = 2,3$. In the forthcoming
sections, for
$q=3$ we provide  examples where $ K^2 = 12$ ( of course the Albanese image is a
surface, as we already observed), while for $q=2$ we show that there are
examples with $K^2 = 18$ (this is not the maximum allowed by the B-M-Y
inequality, yielding $ K^2 \leq 9 \chi = 27 $, but still not so low).

For the case where the degree of the canonical map is $2$ , we recall

\begin{teo} (Debarre's theorem 4.8 [1982] )
Assume that $ 1 \leq q \leq 3$ and that the canonical map has degree $2$: then
$ K^2 \geq 2 p_g + q - 1$.

\end{teo}

\noindent
Finally, we recall the following inequality for fibred surfaces (cf.  Xiao [1987], Konno [1993])
\begin{itemize}
\item
Xiao-Konno inequality: if $f: S \rightarrow B$ is a
fibration to a curve
$B$ of genus $b$, with fibres of genus $g$, 
and without constant moduli, then 
 the slope  $\lambda (f) = (K_S - f^*(K_B))^2  
/ \deg (f_*
\omega_{S|B}) = \frac{K^2_S - 8 (b-1)(g-1)}{\chi(S)- (b-1)(g-1) }$   satisfies
$4 (g-1)/g   \leq \lambda (f) \leq 12 $,  the first inequality being
an inequality iff all
the fibres are hyperelliptic.
\end{itemize}

It follows that $ K^2_S  \leq 12 \chi(S)- 4 (b-1)(g-1) $,
and this  implies that if the Albanese image is a curve of genus $q \geq 2$,
then $ K^2_S  \leq 12 p_g - 12 (q-1) -4 (q-1)(g-1)= 12 p_g - 4 (q-1) ( g +2).
$ 

Therefore, if $p_g = 4$ and if the Albanese image is a curve of genus $q\geq
2$, we obtain indeed $ 16 \leq K^2_S  \leq 48 - 4
(q-1) ( g +2) .$  For $q=3 $ this confirms that we must have $K^2 = 16$ and $g=2$,
a case that we have  already illustrated.  While, if $q=2$ we have 
$  12 \leq K^2_S  \leq 48 - 4 (g +2) $: whence, $g \leq 5$ which is weaker than the inequality
$ g\leq 4 $ we have already obtained.

\section{ Irregular surfaces with $p_g= 4, q = 1$.}

\noindent
Horikawa [1981] proved that for an irregular surface with $ K^2 < 3
\chi $ ($3 \chi=
12$  here),  the Albanese map has a curve as image, and the fibres
have either  genus 2 or they have genus 3 and are hyperelliptic: moreover
only the first possibility occurs  if  $ K^2 < 8/3 \chi .$

\noindent
In this section we shall restrict our
attention to the case $p_g= 4, q = 1$, thus $ \chi=
4$ and  $K^2 = 8 $ is the
smallest value for $K^2$, while $K^2 = 12$ is the smallest value for which we
can have a birational canonical map ( in view of the quoted result by
Horikawa ).

\noindent
Let $a : S \rightarrow A$ be the Albanese map  of $S$.
By the above inequality
for the slope, if $K^2 = 12$ then $g \leq 4$, and if $g=4$ all the
fibres of $a$ are
hyperelliptic. But if the Albanese fibres are hyperelliptic, then the
canonical map $\phi$ factors through the hyperelliptic involution $i$ .
We make
therefore the following

\begin{ass}\label{nohypfib} $p_g= 4, q = 1$, and the general Albanese 
fibre is non
hyperelliptic of genus $g = 3$.
\end{ass}

\noindent
Under the above assumption, let $\sV$ be the
vector bundle on $A $ defined by
\begin{eqnarray} \sV = a_* \omega_S, \end{eqnarray}
which enjoys the base change
property.
$\sV$ has rank $g=3$, and $h^0(\sV) = p_g =4$.  More generally, we can
consider the vector bundles  $\sV_i = a_*(\omega_S^{ \otimes i}) $,
which have zero
$\HH^1$-cohomology groups and have degree 

$\deg (\sV_i ) = \chi + 
{\frac {i(i-1)}{2}} \ K^2 = 4 +
{\frac {i(i-1)}{2}} \ K^2 $.

\noindent
Under the assumption \ref{nohypfib},  we have an exact sequence
\begin{eqnarray} 0 \rightarrow Sym^2 (\sV) \rightarrow \sV_2 
\rightarrow \mathcal C
\rightarrow 0, \end{eqnarray}
where $\mathcal C $ is
a torsion sheaf .

\noindent
Since the degree of $Sym^2 (\sV)$ equals $16 $, we obtain

\begin{prop}\label{degree estimate}  Under the assumption 
\ref{nohypfib}  we have $ K^2 \geq
12 $, equality holding iff there is no hyperelliptic fibre, that is,
$Sym^2(\sV)\cong \sV_2 $.
\end{prop}

\noindent
By the theorem of Fujita [1978] $\sV$ is semipositive, moreover by a corollary
of a theorem of Simpson [1993] observed for instance in 2.1.7 of Zucconi [1994], there is in
general a splitting
\begin{eqnarray} \sV =(\bigoplus _{ \tau \in \Pic(A)_{tors} - \{0\}} L_{\tau} )\oplus
(\bigoplus_i W_i )
\end{eqnarray}
where, as indicated, $L_{\tau} $ is a non trivial
torsion line
bundle, and instead  $W_i$ is an indecomposable bundle of strictly positive
degree.

We are interested mostly in the case where the canonical map $\phi$ is birational;
 since $\sV$ has rank 3, the hypothesis that $\phi$ be birational implies that $\sV$ must be
generically generated by global sections. Thus we make the following

\begin{ass} $\sV$ is
generically generated by global sections . \end{ass}

In particular, there are no
summands of
type $L_{\tau} $ in (5.3) .

Moreover, by Atiyah [1957], Lemma 15, page 430, setting
$\deg(W_i)=d_i$ and $\rank(W_i)=r_i $,
we have $h^0(W_i )= \deg(W_i)=d_i $, and
therefore $d_i \geq r_i$ . Conversely, by loc. cit. Theorem 6, page 433)
and by induction follows

\begin{prop} If $W$ is an indecomposable vector bundle on an elliptic 
curve  of degree $d \geq r= \rank (W)$,
then $W$ is generically generated by global sections.
\end{prop}

Therefore, since $\Sigma_i\, d_i =
4,  \Sigma_i \, r_i = 3$, and the  $r_i$'s, $d_i$ 's  are $>0$, we
have only the
following possibilities for the pairs  $(r_i,d_i)$  which are ordered by the
slope $d/r$:
\begin{itemize}
\item[(i)] (3,4)
\item[(ii)] (2,3) (1,1)
\item[(iii)] (1,2) (2,2)
\item[(iv)] (1,2) (1,1) (1,1).
\end{itemize}

The structure of these bundles is then clear by the
quoted results of Atiyah: for each line bundle $\hol_A(p)$ of degree one,
where $p \in A = \Pic^1(A)$, there are a point  $u \in A $ and line bundles
$ \mathcal L, \mathcal L'
\in  \Pic^0(A)$ such that  $\sV '= \sV \otimes \hol_A(- p) $ is 
respectively equal
to
\begin{eqnarray}\label{cases}
\hbox{(i)} & E_3(u) & \cr
\hbox{(ii)} & E_2(u) \oplus   \mathcal L  & \cr
\hbox{(iii)} & \hol_A(u) \oplus (\mathcal L  \otimes F_2) & \cr
\hbox{(iv)} & \hol_A(u) \oplus \mathcal L  \oplus \mathcal L', & 
\end{eqnarray}

where $ E_1(u) := \hol_A(u)$  and $ E_i(u)$   is defined inductively as
the unique non trivial extension
$$
0 \rightarrow \hol_A \rightarrow E_i(u) \rightarrow E_{i-1}(u)\rightarrow 0,
$$
while $ F_1 := \hol_A$  and
$F_i $ is defined inductively as the unique non trivial extension
$$
0 \rightarrow \hol_A \rightarrow F_i \rightarrow F_{i-1} \rightarrow 0 .
$$
\begin{lem}\label{relCanModel}
         Assume that X is the canonical model of a surface with $ p_g 
=4, q=1, K^2=
12$, and non hyperelliptic Albanese fibres.

Then the relative canonical
map $ \omega : X \rightarrow Proj(\sV) = \PP $  is an embedding. 
Moreover, there
is a  point $p \in A$ such that, setting    $\sV '= \sV \otimes \hol_A(- p) $,
$(\det\sV') \cong  \hol_A(p)$ and $X$ is a divisor in the linear system
$| 4D|$, $D$
being the tautological divisor of $Proj(\sV')$. ( Notice that $p$ is 
defined only
up to 4-torsion) .

Conversely, if $\sV$ is as in (\ref{cases}), any divisor $X$ in $|4D|$ with at
most R.D.P.'s as singularities is the
canonical model of a surface with $ p_g =4, q=1, K^2= 12$, and non
hyperelliptic Albanese
fibres.
\end{lem}

\proof
By \ref{degree estimate}, the relative canonical map is an embedding of the
canonical model $X$ of $S$ if and only if the relative bicanonical map is an
embedding. Let therefore $F$ be a fibre of the Albanese map $a : X
\rightarrow A$.
Since $\sV_2$  enjoys the base change property, we are just asking
whether  $\hol_F ( 2K_X) $
is very ample on F.

By Catanese and Franciosi [1996] or Catanese, Franciosi, Hulek and Reid [1999]
we get that very ampleness
holds provided that
$F$ is 2-connected, i.e., there is no decomposition $F = F_1 + F_2$
with $F_1, F_2$
effective and  with $F_1  F_2 \leq 1$.

If such a decomposition would exist, we claim that we may then assume $F_1  F_2 = 1$.

 Otherwise $F_1  F_2 \leq 0$ and, 
since also $ F_i^2 \leq 0$, it follows by Zariski's Lemma that $ F = 2 F_1$.
Since $F_1$ has genus $2$ the exact sequence 
$$ \hol_{F_1} ( K_X - F_1) \ra  \hol_{F} ( K_X ) \ra \hol_{F_1} ( K_X) \ra 0$$
shows easily that the hyperelliptic involution on $F_1$ extends toa hyperelliptic involution on
$F$, a contradiction. 

Since the genus of $F$ equals $3$, $K_X F = 4 = K_X F_1 + K_X F_2 $.
Since moreover
$ K_X F_i \equiv F_i ^2 \ (mod \ 2) \equiv - F_1  F_2  \ (mod \ 2)$, we
may assume w.l.o.g.
that $ K_X F_1  = 1$. So $F_1$ is an elliptic tail, while $F_2$ has genus $2$.

 More precisely, we
have
$h^0(\hol_{F_1}(K_X)) = 1$, $h^0(\hol_{F_1}(2 K_X)) = 2$,
contradicting the fact that 
$ Sym^2 (H^0(K_F)) = (H^0( 2 K_F))$, since in fact 
$h^1(\hol_{F_2}(2 K_X - F_1)) = 0.$

We proved now that $X$ is embedded in $\PP$, whence it follows that the surjection $
Sym^4 (\sV) \rightarrow \sV_4$
has as kernel an invertible sheaf $L$ on the elliptic curve $A$.
The exact sequence
$$
0 \rightarrow L \rightarrow Sym^4 (\sV) \rightarrow \sV_4 \rightarrow 0,
$$
and the easy calculation: $ \deg (Sym^4 (\sV)) = 80, \deg(\sV_4) = 
76$ shows that
$ \deg(L) = 4 $, so there is point $p \in A$ such that $ L = \hol_A (4p)$,
and we have the following linear equivalence in $\PP$: $ X \equiv 4H
- 4 F$, where
$H$ is the tautological hyperplane divisor, and $F$ is the fibre  of
$a$ over $p$.

If we write $\PP$ as $ Proj (\sV') = Proj (\sV(-p))$, and let $D$ be the
corresponding
   hyperplane divisor, then $X$ is a divisor in the linear system $|4 D|$.

Conversely, since $ K_{\PP} \equiv -3 H + \omega^* (\det \sV) $, if we
choose a divisor
$X \equiv  4 H - \omega^* (\det \sV) $ we get that $ K_X \equiv H$, so
that $\sV$ is the
direct image of the canonical sheaf $ \hol_X(K_X)$.
It is then clear that $L =\omega^* (\det \sV) $, and since we chose 
$p$ such that
$ X \equiv 4H - 4 F$, $ K_X \equiv  D + \omega^* (\det \sV') \equiv D + F$,
   we have proven that $ (\det \sV')  \cong \hol_A(p)$.

Moreover, we have that $p_g(X) = h^0(\sV) = 4 $, $ q = 1  + h^1(\sV) 
= 1 + 0 = 1$,
while $ K_X ^2 = (D + F)^2 (4 D) $. By the Leray-Hirsch formula we get
$ D^2 = FD$, moreover $ D^2 F = 1 $, whence $ K_X ^2 = 4 ( 1 +2) = 12$.
\QED

\begin{rmk}
Observe that $ det(\sV ) \cong \hol_A (4p) $, whence in the  notation of the
previous lemma we have $ det(\sV') \cong \hol_A (p) $. Thus 
$ p \equiv u$ in case i), while $ p \equiv u + \mathcal L$,
$ p \equiv u + 2 \mathcal L$, $ p \equiv u + \mathcal L  + \mathcal L'$,
in the respective cases ii), iii), iv).

It follows that the pair $(A, \sV)$ has $1$ modulus  in case i), $2$ moduli
in case ii), while we are going to show next that the pair has $1$ modulus in
case iii), and $1$ or $2$ moduli in the last case iv).
\end{rmk}

There remains as a first problem the question about the existence of the surfaces
under consideration, that is, whether the general element in the linear system
$|4D|$ has only Rational Double Points as singularities.
The result is a consequence of techniques developed in 
Catanese and Ciliberto [1993].

\begin{prop}\label{RDPs}
Let $\sV$ be as in (\ref{cases}) a rank $3$ bundle over an elliptic curve, and $X$ as
in Lemma \ref{relCanModel} a general divisor in the linear system $|4D|$ on $\PP =
Proj(\sV)$. Then $X$ is smooth in cases i) and ii).
Instead, in case iii), the general element $X$ has only Rational Double
Points as singularities  if and only if $\mathcal L^4 \cong \hol_A$.

Finally, in case iv), the general element $X$ has  Rational Double
Points as singularities  if and only if one of the bundles 
$\mathcal L^k \otimes {\mathcal L' }^{4-k}$ is trivial.
\end{prop}

\proof

In case i), the linear system $|4D|$ is very ample on $\PP$ by Theorem 1.21
of Catanese and Ciliberto [1993], so a general $X$ is smooth by Bertini's Theorem.

Notice that to apply Bertini's theorem it suffices to show that the
general element of the linear system
$|4D|$ is smooth along the base locus of $|4D|$. To show that $|4D|$ is
base point free  is in turn sufficient the condition that the vector
bundle $Sym^4(\sV')$ be generated by global sections.

In case ii), we observe that $Sym^4(\sV')$ is a direct sum
$$ \bigoplus_{k=1,..4} Sym^k(E_2(v_k)) \bigoplus \mathcal L^4,$$ where the
$v_k$'s are suitable points on $A$. The
  symmetric powers with $ k \geq 2$ are generated by global sections by virtue
of Theorem 1.18 of Catanese and Ciliberto [1993].

Whereas a bundle of the form $E_2(v_k)$ has only one section, which is
nowhere vanishing. Therefore,  the base locus of $|4D|$ is contained in the
section $\Delta$ of $\PP$ dual to the subbundle $E_2(u)$, and there, if $x,y$
are local equations for $\Delta$, then the Taylor development of the equation
of a divisor in $|4D|$ only fails at most to have a term of type $x$.
Therefore in this case the general element $X$ is smooth.

Let us then consider case iv): it is immediate to remark
that the linear system $|4D|$ has as fixed part the projective $\PP^1$-subbundle
$\PP '$ annihilated by $\hol_A(u)$ in the case where no
line bundle $\mathcal L^k \otimes {\mathcal L' }^{4-k}$ is trivial.

This cannot occur, so assume that one of such line bundles is trivial.

{\bf Case iv) -(I''):} {\it $ {\mathcal L'}^{4}$ and $\mathcal L^{4}$ are trivial.}

This is the easy case where   $|4D|$ has no base points, whence Bertini's
theorem applies.

{\bf Case iv) -(I'):} {\it ${\mathcal L'}^{4}$ is not trivial but $\mathcal L^{4}$
is trivial.}

This is the case where the base locus of $|4D|$ is the section $\Delta$
annihilitated by $\hol_A(u) \oplus \mathcal L$.

At each point of $\Delta$ exists then a term of order $1$ in the Taylor
expansion of the equation $f$ of $X$, except at the point $t'=0$, where
$t'=0$ lies over the unique point $v \in A$ such that $ v \equiv u + 3
\mathcal L'$.

In this point we have then local coordinates $x,z, t'$ with $\Delta = \{ z = x =0\}$ and we
surely get monomials $z^2, z t', z x$ ( $zx$ corresponds to the fact that the unique
section of $\hol_A(u) \otimes \mathcal L \otimes{\mathcal L'}^2$ does not
vanish in our point $t'=0$, otherwise $\mathcal L \cong {\mathcal L'}$,
contradicting our assumption (iv, I')). Since moreover we get the monomial $x^4$,
we certainly obtain for general $X$ at worst a Rational Double Point of type
$A_3$.

{\bf Case iv) -(II):} {\it ${\mathcal L'}^{4}$ and $\mathcal L^{4}$ are not
trivial but $\mathcal L^2 \otimes {\mathcal L' }^2$ is  trivial
.}

In this case the base locus of
$|4D|$ is given by the two sections $\Delta$, $\Delta '$, where 
$\Delta '$ is annihilated by $\hol_A(u) \oplus \mathcal L'.$

In this case, by symmetry, let us study the singularity of a general $X$ along $\Delta$.

We get, as in the previous case, a section $ z t'$, and sections $z^2, x^2$ , thus a 
singularity of type $A_1$ at worst.

{\bf Case iv) -(III):}  {\it ${\mathcal L'}^{4}$, 
$\mathcal L^2 \otimes {\mathcal L' }^2$ and
$\mathcal L^{4}$ are not trivial but $\mathcal L^3 \otimes {\mathcal L' }$ is trivial.}

In this case again the base locus of
$|4D|$ is given by the two sections $\Delta$, $\Delta '$.
For $\Delta$ we get monomials of type $ x^3, z^2, z t'$, thus a singularity of type $A_2$ at
worst, while for $\Delta'$ we get monomials of type $ x, z^2, z t'$, thus no singularity at all
for general $X$ along $\Delta '$.

We can finally analyse case iii).

To this purpose, recall 

\begin{teo}(Atiyah's Theorem $9$ in [57])
Let $F_2 $ be the indecomposable bundle on an elliptic curve with trivial
determinant and of rank $2$: then $Sym^k(F_2) \cong F_{k+1} $.
\end{teo}

We observe then that the tensor product of $F_r$ with a line bundle $\mathcal
M$ is generated by  global sections if $ \deg (\mathcal M) \geq 2$, whereas,
if $ d:= \deg (\mathcal M) \leq 1$ , $F_r \otimes \mathcal M$ is generated by
global sections outside a unique point where all sections vanish if $d=1$,
whereas for $d=0$,  $F_r \otimes \mathcal M$ has no global sections unless
the line bundle  $ \mathcal M$ is trivial. In this last case, there is only
one non zero section, which vanishes nowhere.

After these remarks, it is clear that, in case iii), if $\mathcal L^4$ is non
trivial, then the fixed part of $|4D|$ contains the $\PP^1$-bundle $\PP '$
annihilated by $\hol_A(u)$. Otherwise, the base locus consists of a section
$\Delta$ and a fibre $F_v$ of $\PP '$.

At the intersection point of these two curves, we have local coordinates
$z, x,t$ such that $z=0$ defines $\PP '$, $z=x=0$ defines $\Delta$,
$z=t=0$ defines $F_v$.

In  the Taylor expansion of the equation of $X$ we get $x^4$, $z t$, $z^2$,
therefore, by an argument we already used, in this point we get at worst
a singularity of type $A_3$, while the other points of the base locus are
smooth for general $X$.

\QED

We derive from the previous result some preliminary information on the moduli 
space of the above surfaces.

We need however to slightly simplify our previous presentation.
Observe therefore that we can exchange the roles of ${ \mathcal L} $, ${ \mathcal L'} $, and we
can tensor $\sV '$ by a line bundle of $4$-torsion.

Therefore, in the last case iv), we may simplify our treatment to consider only the subcases

\begin{itemize}
\item
(iv, I) ${ \mathcal L} \cong \hol_A$
\item
(iv, II) ${ \mathcal L' } \cong { \mathcal L}^{-1}$, ${ \mathcal L}$ not of 4-torsion
\item
(iv, III) ${ \mathcal L' } \cong { \mathcal L}^{-3}$, ${ \mathcal L}$ not of 4-torsion
\end{itemize}
( the reader may in fact observe that cases (iv, I') and (iv, I'') are just
special subcases of (iv, I)).

\begin{cor}
Consider the open set $\mathcal M$ of the moduli space of the surfaces with
$p_g=4, q=1, K^2 = 12$ such that assumptions 5.1 and 5.3 are verified.

Then $\mathcal M$ consists of the following  $10$ locally closed subsets:

\begin{itemize}
\item
$\mathcal M (i)$, of dimension $20$, corresponding to case i)
\item
$\mathcal M (ii,0)$ corresponding to the case $\mathcal L ^4 \cong \hol_A$,
and $\mathcal M (ii,1)$ corresponding to the case $\mathcal L ^4 \ncong \hol_A$,
both of dimension $19$
\item
 $\mathcal M (iii)$ corresponding to case iii), of dimension $18$
\item
 $\mathcal M (iv,I)$ of dimension $18$, corresponding to the case where 
$\mathcal L  \cong \hol_A$, but  ${ \mathcal L'} $ is neither of 3-torsion nor
of 4-torsion
\item
$\mathcal M (iv,II)$ of dimension $19$, corresponding to the case iv), II),
(${ \mathcal L}$ not of 4-torsion)
\item
$\mathcal M (iv,III)$ of dimension $19$, corresponding to the case iv), III), 
(${ \mathcal L}$ not of 4-torsion)
\item
$\mathcal M (iv,I,1/4)$ of dimension $18$, corresponding to case iv),I), where
we may assume $\mathcal L  \cong \hol_A$ and $\mathcal L ' $ of $4$-torsion
but not of $2$-torsion
\item
$\mathcal M (iv,I,1/2)$ of dimension $19$, corresponding to case iv),I), where
we may assume $\mathcal L  \cong \hol_A$ and $\mathcal L ' $ of $2$-torsion,
but non trivial
\item
$\mathcal M (iv,I,1/3)$ of dimension $18$, corresponding to case iv),I), where
we may assume $\mathcal L  \cong \hol_A$ and $\mathcal L ' $ of $3$-torsion
and non trivial
\item
$\mathcal M (iv,I,1)$ of dimension $19$, corresponding to case iv),I), where
we may assume $\mathcal L  \cong \hol_A$ and also $\mathcal L ' \cong \hol_A. $ 
\end{itemize}

\end{cor}

\Proof
Observe that the moduli space of elliptic curves together with a torsion sheaf
of torsion precisely $n$ is irreducible of dimension $1$.

Observe moreover that the hypersurface $X$ ( the canonical model of $S$) moves
in a linear system $|4D|$ in $\PP$, whose dimension is given by
$h^0(Sym^4(\sV ')) -1$. By Riemann-Roch,  since $\sV'$ has rank $=3$ and
degree $=1$,
 we have that
$h^0(Sym^4(\sV '))= 20 + h^1(Sym^4(\sV '))$. 

Moreover we have that the dimension of each stratum, since each surface of general type has a
finite automorphism group, equals $1 +h^0(Sym^4(\sV ')) -h^0(End(\sV ')) $. 
 This justifies the assertion about the
dimensions.

Furthermore, we may observe that the conditions that a vector bundle be
indecomposable is an open one, while the condition that a line bundle be
of $n$-torsion is a closed one.

\QED

The previous corollary allows us to conclude that our moduli space is irreducible:
we use for this purpose a lower bound for the dimension of the moduli space which is a
consequence of a general principle stated by Ziv Ran [1995], and turned into a precise theorem
by Herb Clemens [2000].

\begin{teo}
The open set of the moduli space of surfaces with $ p_g
=4, q=1, K^2=
12$,  with non hyperelliptic Albanese fibres,
   and  with $\sV$ as in (\ref{cases}), i.e.,  generically generated by
global sections, is irreducible of dimension $20$. 

\end{teo} 
\proof
In view of  the previous corollary, our moduli space has a stratification by locally closed
sets, of which one only, $\mathcal M (i)$, has dimension $20$, while the others have strictly
smaller dimension. Since $\mathcal M (i)$ is clearly irreducible, it suffices to show that the
dimension of the moduli space is at least 20 in each point.
Equivalently, since the germ of the moduli space at the point corresponding to the surface $S$ is
the quotient of the  base of the Kuranishi  family  of $S$ by the finite group of automorphisms
of $S$, it suffices to show that the dimension of the Kuranishi family is at least $20$.  

Now, in any case, the dimension of the Kuranishi family is always at least
$$ h^1 (S,T_S) - {\rm dim (Obs}(S)) ,$$ but in this case the obstruction space $Obs(S)$ is not
the full cohomology group $ H^2 (S,T_S)$. Because we have a natural Hodge bilinear map 
$$ \gamma : H^0 ( S, \Omega^1_S) \times H^0 ( S, \Omega^2_S) \rightarrow H^0 ( S,\Omega^1_S
\otimes
\Omega^2_S),$$ and the natural subspace $ H := {\rm Im }(H^0 ( S, \Omega^1_S) \otimes H^0 ( S,
\Omega^2_S)) \subset H^0 ( S,\Omega^1_S \otimes \Omega^2_S)$ determines  by Serre duality a
quotient map $  \gamma^{\vee}:H^2 (S,T_S) \rightarrow H^{\vee}$. By Theorem 10.1 of Clemens
[2000] we have that   $ \gamma^{\vee} (Obs(S)) = 0.$ 

Since in this case it is obvious that ${\rm dim} (H) = 4$, $\gamma$ being non degenerate, it
follows that the base of the Kuranishi family has dimension $\geq  - \chi (S,T_S) + 4 = 10
\chi(S) - 2 K_S^2 + 4 = 20.$
\QED

\begin{teo}
         Assume that X is the canonical model of a surface with $ p_g
=4, q=1, K^2=
12$, and non hyperelliptic Albanese fibres,
   and $\sV$ is as in (\ref{cases}), i.e.,  generically generated by
global sections . Then in cases (i), (ii) the canonical map $\phi$ is
always a birational morphism, whereas in the other cases  $\phi$  is
birational for a general choice of $X$ in the given linear system.
The case $\deg (\varphi) = 3$ never occurs.
\end{teo}
\proof
Since $\sV$ is   generically generated by
global sections, and the general fibre $F_a$ is a non hyperelliptic curve of
genus
$3$, it follows that $F_a$ maps isomorphically to a plane quartic curve
$\Gamma_a$.
Let $H_a$ be  the plane containing $\Gamma_a$: since $K_X \equiv D + F$,
then the pull-back divisor of $H_a$ splits as $F_a + D_{-a}$, where $ D_{-a}
\sim D$.

Since $12 \geq \deg \phi \cdot  \deg \Sigma$, and there are plane sections
$H_a$ which intersect $\Sigma$ in a curve containing $\Gamma_a$, $\deg
\Sigma \geq 4$,  and the only possibility to exclude is that $\deg \phi = 2$ or
$3$.

In the case where  $\deg \phi = 2$ let $\imath : X \ra X$ be the corresponding
biregular involution.

$\imath$ acts also on the Albanese variety $A$, and in a non trivial way,
since a general fibre $F_a$ is embedded by the canonical map $\phi$, and let
us then denote $a':= \imath (a)$.

If $\imath$ had no fixpoints on $A$, then $X \rightarrow X / \imath := Y$ 
would be unramified, so that $K^2_Y = 6$, $\chi (\hol_Y )= 2$, $q(Y) = 1$,
whence $p_g (Y) = 2$, contradicting the fact that $\phi$ factors through
$Y$. We may therefore assume that $ a' = -a$ for a suitable choice of the origin in $A$.

Therefore, the inverse image of $H_a$ contains
$F_a + F_{a'}$, and we can write a linear equivalence $ K_X \equiv F_a +
F_{a'} + C_a$, where
$C_a$ is effective and $C_a = C_{-a}$.

Observe that $K_X \cdot F = 4$, $K_X \cdot C_a
= 4$, $C_a \cdot F =
4$ whence  $C_a$ is not vertical for the Albanese map. Moreover, 
$ 12 = K_X^2 = ( 2 F + C_a)^2= 16 + C_a^2$, whence $C_a^2 = -4$.

In particular the algebraic system $C_a$ has a fixed part. 

We obtain a contradiction as follows.

First of all, since $| K_X - F_a -
F_{a'} | \neq \emptyset $, we get $H^0(A, \sV (-a-a')) \neq 0$,
and this leaves out only the cases (\ref{cases}) iii) and iv), and moreover
with $ a+ a' \equiv u + p$ on $A$.

We saw that $A/ \imath \cong \PP^1$, so that all the
curves $C_a$ are linearly equivalent. Indeed, a closer look reveals that all
the curves $C_a$ are the intersection of $X$ with a fixed $\PP^1$-subbundle
of $\PP$, thus we may consider the curve $ C = C_a , \forall a \in A$.

The curve $C$ maps to a line $L$ under the two dimensional linear system
corresponding to $H^0(A,W)$, where we write $\sV = \hol_A(u +p) \oplus W$.

Before we further investigate the geometry of the situation, remark that
$\imath$ acts equivariantly on $X$ and $A$, therefore $\sV$ is isomorphic
to $\imath ^*(\sV)$ and indeed we have an action of $\imath $ on $\sV$.

This however implies that $\mathcal L$ is of $2$-torsion in case iii), while
in case iv) $\mathcal L \cong  - \mathcal L'$. Once these conditions are
satisfied, it is clear that we have an involution $\imath$ on $\PP$ and that
the system $| \hol_{\PP}(1)|$ is invariant, but it remains to be seen whether
the hypersurface $X$ is also $\imath$-invariant (notice that the involution
is completely determined by the four fixed points $O$ such that $ 2O \equiv u
+ p$).

It is easy to verify that for a general choice of $X$ in $|4D|$, this does
not hold.

{\bf CLAIM: $\deg (\phi) = 3$  NEVER OCCURS.}

Consider in fact the possibility that $deg (\phi) = 3$: then $\Gamma_a$ is a
full hyperplane section of $\Sigma$, and $K_X$ is base-point free
(in general,  if $|K_X| = |M| + \Psi$, with $\Psi$ a non trivial fixed part,
then $M^2 =  K_X^2 - K_X \cdot \Psi - M \cdot \Psi <  K_X^2$, if then $|M|$
has base points, then $M^2 > \deg (\phi) \cdot \deg(\Sigma)$: while here $
K_X^2 = 12 = \deg (\phi) \cdot \deg(\Sigma)$).

Observe that the surface $\Sigma$ is normal, since it has a smooth hyperplane
section.

Let $\pi: \tilde \Sigma \ra \Sigma$ be a minimal resolution of $\Sigma$ and
denote by
$\tilde X$ a minimal resolution of the fibre product $\tilde \Sigma
\times_{\Sigma} X$: since $\tilde X$ is birational to $X$, $R^1 p_*
\hol_{\tilde X} =0$, $p:
\tilde X \ra \Sigma$ being the composite morphism. Whence it follows that
  $R^1 \pi _* \hol_{\tilde \Sigma} =0$, i.e., $\Sigma$ has only rational
double points as singularities, and  $\tilde \Sigma$ is a smooth $K3$-surface.

We will now consider the ramification formula for $\phi$.  Let $B$ be the
reduced  branch divisor of $\phi$,
set $\phi^* (B) = R + R'$, $R$ being the ramification divisor, and observe
that $R' \geq 1/2 R$.
The fact that $\Sigma$ is a $K3$ with R.D.P.'s implies that $R \equiv K_X$,
i.e., there is a hyperplane divisor $H$ with $R = \phi^* (H)$.

On the other hand, since $ \deg (\phi) =3$ it follows that $\phi_*(R_{red})=B$,
whence $B$ is  the reduced divisor of the plane section which pulls back to
$H$: this is a contradiction  since then $R = \phi^* (H)  \geq \phi^*
(B) = R + R'$, while
$R' > 0$ ( this follows since $R > 0$, and $R' \geq 1/2 R$).
\hfill Q.E.D.

\begin{rmk}

We saw that we have several strata of the  above irreducible moduli space. 
The stratum of maximal dimension, such that the moduli space is just its closure
will be called the  'Main Stratum',  and we shall say that the surfaces
which belong to this Main Stratum are of the "Main Stream".

It is certainly, as we shall see , the one which is most interesting and related
to the geometry of elliptic space curves of degree $4$.

We should also remark that a detailed and more general study of surfaces with
irregularity $ q=1$ and with $ K^2 = 3 \chi$ was undertaken in the $1996$
Thesis of T. Takahashi. However his results are weaker than ours in the case
where $p_g = 4$, so we could not use this reference. 
\end{rmk} 

We finally come to a discussion of the geometry of the surfaces of the "Main Stream"
( case i) ).

Let $A$ be an elliptic curve of degree $4$ in $\PP^3$. Then, as it is well
known, $A$ is the complete intersection of $2$ quadric surfaces $Q, Q'$.

We may indeed without loss of generality assume that the pencil of quadrics
be Heisenberg invariant, in other words that:

$$ A = \{ (x) | x_0^2 + x_2^2 - \lambda^2 x_1 x_3 =0, x_1^2 + x_3^2 - \lambda^2 x_0 x_2 =0   
\}. $$

It is also well known (Atiyah [1957], Catanese and Ciliberto [1988]) that in case
i) the projective bundle $\PP$ is nothing else than the triple symmetric product of
the elliptic curve $A$, $\PP = A^{(3)}$.

In this context the canonical mapping of $X$ is induced by a morphism 
$\phi: \PP \ra (\PP^3)^{\vee} $ which can be explained without formulae as
follows: consider a point $P$ of $\PP$, i.e., $P$ is a divisor of degree $3$ on
$A$. Then there is a unique plane $\phi (P)$ containing this divisor.

This geometric explanation shows that the degree of $\phi$ is $4$ ( as it had to
be, since $F$ being a fibre of the Albanese map $a$, $(D+F)^3 = 4$ by the
Leray-Hirsch formula which we already mentioned).

In fact the reason of the above is that 
\begin{itemize}
\item
the projection onto the elliptic curve
(the Albanese map) associates to a divisor $P = P_1 + P_2 + P_3$ the sum of the
three points $P_1 , P_2 , P_3$ in the elliptic curve $A$,
\item
  the tautological
divisor $D$ on $\PP$ consists of the divisors where $P_1$ is fixed (whence,
$D^3 = 1$). 
\item
Let $I= {\rm Proj} (\sT_{\PP^3})$ be the incidence correspondence, $I \subset \PP^3
\times(\PP^3)^{\vee}$: then we claim that  $I \cap ( A
\times(\PP^3)^{\vee}) = A^{(3)}$.
\proof
The isomorphism is given by $ (x, h) \ra  {\rm div}_A(h) - x . $
\qed
\bigskip

Observe only that the first projection does not correspond precisely to the
Albanese map, but only to the composition of the Albanese map with an involution
of $A$, since to a divisor $P = P_1 + P_2 + P_3$  corresponds the point $x$ such
that $x + P_1 + P_2 + P_3$ is linearly equivalent to the hyperplane divisor of
$A$.
\end{itemize}
\begin{itemize}
\item
We claim that the second projection is given by the linear system $|D+F|$.
\proof
Any hyperplane in $(\PP^3)^{\vee}$ is the hyperplane $H_x$ dual to a point $x\in
\PP^3$. Let $x \in A$: then the inverse image of $H_x$ is given by the divisors
$P'$ of degree $3$ on $A$ such that  $P'_1 , P'_2 , P'_3$ span a plane containing
$x$. Thus, we have two possibilities: either we take the divisors $P$ such that
$P + x$ is linearly equivalent to the hyperplane divisor on $A$, and thus we get
a fibre $F$, or we simply take the divisors $P'$ of degree $3$ for which $ P' \geq
x$, i.e., we get a divisor of type $D$.
\end{itemize}

Set for convenience $W : = \PP = A^{(3)} $ and observe that the pull back $H_2$
of the  hyperplane divisor  in $(\PP^3)^{\vee}$ is thus linearly equivalent to
$D+F$. Observe also that the pull back of the hyperplane divisor  in
$(\PP^3)$ is linearly equivalent to $4F$. Therefore, the desired canonical model
$ X \subset W$ is in the linear system $ |4 D | = | 4 H_2 - H_1  |$.

We can perhaps summarize these observations as follows:

\begin{prop}
The canonical model of a surface with $p_g=4, q=1, K^2 = 12$ of the Main Stream,
i.e., of type i), is a divisor of bidegree $(-1, 4)$ on the variety
 $W$  given by the intersection of the incidence variety
$I \subset \PP^3
\times(\PP^3)^{\vee}$ ( itself a divisor of bidegree $(1,1)$) with the pull
back of the elliptic curve $A$ under the first projection.

Thus $W$ is a complete intersection of type $(1,1), (2,0)(2,0)$, but $X$ is not
a complete intersection. The canonical divisor on $X$ is induced by the divisor
of bidegree $(0,1)$ on  $ \PP^3
\times(\PP^3)^{\vee}$.
\end{prop}

From this it is easy to produce equations of explicit examples of these surfaces
via computer algebra. The method is based on the following
 
\begin{rmk}
Since $W$ is a complete intersection in $ M:= \PP^3 \times(\PP^3)^{\vee}$ it follows
easily that the restriction homomorphism $\HH^0 ( \hol_M (n,4)) \ra  \HH^0 (
\hol_W (n,4))$ is surjective as soon as $n \geq 2$.

Fix therefore a divisor $\bar B$ in $|\hol_M(3,0)|$, for instance the pull back of
the three planes $ \{ x| x_0 x_1 x_2 =0 \}.$ Then the linear system $|4D|$ on $W$
is the residual system $|\HH^0 (\hol_W (3,4) (- \bar B))|.$
\end{rmk}

\bigskip

\section{Surfaces with $p_g=4, q=3$ and canonical map of degree $1$ or $2$}

In this section we shall consider surfaces with $p_g=4, q=3, K^2= 12$,
contained in an Abelian 3-fold as a polarization of type $(1,1,2)$:
we will first show that this family is stable by small deformation.

Later, we will show that for a general such surface the canonical map is a
birational morphism onto a surface of degree twelve in $\PP^3$,
whereas, for all
the surfaces which are the pull back of a theta divisor on a principally
polarized Abelian 3-fold, then the canonical map is of degree
$2$ onto an interesting sextic surface.

More precisely, our situation will be as follows: we let $J$ be a
principally polarized Abelian variety of dimension $3$, which is the
Jacobian of a curve $C$ of genus $3$, and we let $\Theta$ be its principal
polarization. We let $\pi: A  \rightarrow J$  be an isogeny of degree $2$
and $S$ a smooth divisor in the complete linear system $|\pi^* \Theta|$
associated to the pull back $\pi^* \Theta$.

Since at some step we will also need theta
functions, we represent the Jacobian variety $J$ as $ J = \C^3 / \Z^3 + \Omega
\Z^3 : = \C^3 / \Lambda (\Omega)$, with $\Omega$ in the Siegel upper
half-space (we
have thus already represented the theta divisor as a symmetric divisor with
respect to the origin in $J$).

We set then $ A = \C^3 / \Lambda $ with
$ \Lambda  \subset \Lambda (\Omega)$  of index $2$ dual under the
symplectic pairing
to $ \Lambda (\Omega) + \Z b$,  with $ b \in 1/2 \Z^3$ ( e.g., $ b = 1/2 e_3$).
\bigskip
\begin{prop}\label{even functions}
A basis of $ \HH^0 ( A, \hol_A( \pi^* \Theta ))$ is given by even functions.
\end{prop}
\proof
Let $c \in \{ 0, b \}$, and consider the  basis given by the following two
elements:
$$\theta [0,c] ( z, \Omega) : = \sum_{n \in
\Z^3} exp ( 2 \pi i \ (1/2 \ ^t n \Omega n + ^t n ( z + c)    )).$$
An elementary calculation shows that
$$\begin{matrix}
\theta [0,c] (-  z, \Omega) & = &
\sum_{n \in
\Z^3} exp ( 2 \pi i \ (1/2 \ ^t n \Omega n + ^t n ( - z + c)    ))\cr
& = &\sum_{m = -n \in
\Z^3} exp ( 2 \pi i \ (1/2 \  ^t m \Omega m + ^t m ( z - c)    )) \cr &=&
\theta
[0,c] ( z, \Omega)   \cr
\end{matrix}$$
since,  $ \forall m \in \Z^3$, $exp ( 2
\pi i ( ^t m (-2c)) ) = 1 $.

\QED

\begin{prop}\label{numData}
Let $S$ be a smooth surface in a polarization of type $(1,1,2)$ in an Abelian
3-fold $A$. Then the invariants of $S$ are  $p_g=4, q=3, K^2= 12$.
\end{prop}
\Proof
Let us consider the exact sequence
$$0 \ra \hol_A \ra \hol_A(S) \ra \omega_S \ra 0$$
and observe that $\HH^i (\hol_A(S)) = 0 $ for $i=1,2$. Whence, $p_g = h^0(\omega_S) = 4$, and
$ q :=h^1 (\hol_S) $($= h^1 (\omega_S)$ by Serre duality) $=3$.

Moreover, we have $K_S^2 = S^3 =  12$.
\QED
\begin{prop}
Let $S$ be a smooth surface in a polarization of type $(1,1,2)$ in an Abelian
3-fold $A$. Then any small deformation is a surface of the same kind.
\end{prop}
\Proof
Since the canonical divisor of $A$ is trivial, the normal bundle of $S$ in $A$
is $ N_S = \omega_S$, whence its cohomology groups have respective dimensions
$h^0(N_S)= 4,h^1(N_S)= 3,h^2(N_S)= 1. $

The tangent sheaf sequence reads out as follows:
$$0 \ra \sT_S \ra 
\sT_A \otimes \hol_S \cong \hol_S^3 \ra N_S \ra 0,$$
whose exact 
cohomology sequence is:
$$0 \ra \HH^0(N_S)/\HH^0(\hol_S^3) \cong \C 
\ra \HH^1(\sT_S) \ra 
\HH^1(\sT_A \otimes
\hol_S )\ra 
\HH^1(N_S) \ra \HH^2 (\sT_S)
\ra ...$$
We get a smooth 
$7$-dimensional family by varying $A$ in its $6$-dimensional
local 
moduli space (Siegel's upper half space), and $S$ in the 
corresponding
$1$-dimensional linear system.

This family will be 
shown to coincide with the Kuranishi family once 
we prove 
that
$\HH^1(\sT_S)$ has dimension $7$, or, equivalently, (since 
$\HH^1(\sT_A \otimes
\hol_S )\cong \HH^1(\hol_S^3) \cong \C^9$) 
we show the surjectivity of
$\HH^1(\sT_A \otimes
\hol_S )\ra 
\HH^1(N_S)$.

To understand this map, consider an element $ 
\sum_{i=1,2,3} \xi_i 
\otimes \psi_i
\in  \HH^1(\sT_A \otimes 
\hol_S )$, where $\xi_1, \xi_2, \xi_3$ yield a basis
of 
$\HH^0(\sT_A)$, $\psi_i \in  \HH^1(\hol_S) \cong  \HH^1(\hol_A) 
$.

Let $\{ U_{\alpha} \}$ be an open cover of $A$ such that $S \cap 
U_{\alpha} =
div(f_{\alpha})$, and let $f_{\alpha} = g_{\alpha, 
\beta} f_{\beta}$  in
$ U_{\alpha} \cap U_{\beta}$:

then the image 
of  $ \sum_{i=1,2,3} \xi_i \otimes \psi_i$ is given by
$ 
\sum_{i=1,2,3} \xi_i (f_{\alpha}) \otimes \psi_i$.

We use moreover the isomorphism $ \HH^1(N_S) \cong  \HH^2(\hol_A) $: 
since for a
vector field
$\xi$ we have $\xi (f_{\alpha}) = g_{\alpha, \beta} \ \xi (f_{\beta}) ( mod
\  f_{\beta})$, the image of $ \sum_{i=1,2,3} \xi_i \otimes \psi_i$ into
$\HH^2(\hol_A) $ is the cohomology class $ \sum_{i=1,2,3} \xi_i
(g_{\alpha, \beta})
\cup \psi_i$.

We are quickly done, since
\begin{itemize}
\item
the map $ \xi \in \HH^0(\sT_A) \ra \xi (g_{\alpha,
\beta}) \in \HH^1(\hol_A)$ is an isomorphism, being the tangent map at the
origin of the isogeny $\tau: A \ra \Pic(A)$ such that $\tau(x) = S - (S+x)$
\item
$ \HH^1(\hol_A)  \cup \HH^1(\hol_A) \ra \HH^2(\hol_A) $ is onto.
\end{itemize}
\QED

   \begin{teo}
Let $S$ be a smooth divisor yielding a polarization of type $(1,1,2)$
on an Abelian
$3$-fold: then the canonical map of $S$ is  in general a birational
morphism onto a
surface $\Sigma$ of degree $12$.

In the special case where $S$ is the inverse image of the theta divisor in a
principally polarized Abelian 3-fold, the canonical map is a degree
two morphism
onto a sextic surface $\Sigma$ in $\PP^3$.
In this case, the singularities of $\Sigma$ are in general: a plane cubic
$\Gamma$ which is a double curve  of nodal type for $\Sigma$ and,
moreover, a
strictly even set of $32$ nodes for $\Sigma$. Also, in this special case,
the normalization of
$\Sigma$ is in fact the quotient of $S$ by an involution $\iota $ on $A$ having
only isolated fixed points (on $A$), of which exactly $32$ lie on $S$.
   \end{teo}

\Proof
Observe that the natural map $ \HH^0(\Omega^2_A) \ra \HH^0(\Omega^2_S)$
is injective
because $S$ is not a subabelian variety, moreover we get in this way a linear
subsystem  of $|K_S|$ which is base point free, since $S$ embeds into $A$.

It is easy to observe that each translation, and also each involution
$\iota $ with
fixed points on $A$ (multiplication by $-1$ for a suitable choice of an origin)
acts trivially on the vector space $ \HH^0(\Omega^2_A)$.

On the other hand, considering the exact sequence in Prop. \ref{numData},

$$0 \ra \omega_A \ra \omega_A(S) \ra \omega_S \ra 0$$

we see that the $3$-dimensional system generated by $ \HH^0(\Omega^2_A)$ maps
isomorphically to $ \HH^1(\omega_A)$, whereas $\HH^0(\hol_A(S)) \cong
\HH^0(\omega_A(S))$ maps to $\HH^0(\omega_S)$ under the following explicit map
$$ f(z)  \ra f(z) (dz_1 \wedge dz_2 \wedge dz_3 ) / d \theta(z),$$
where $S = div( \theta(z)) $, $f(z),\theta(z) \in  \HH^0(\hol_A(S))$ 

are expressed
by even functions, and  $ (dz_1 \wedge dz_2 \wedge 
dz_3 ) / d 
\theta(z)$ stands for
the Poincare' Residuum $\eta:= 
\eta_i := (dz_1 \wedge dz_2 \wedge 
dz_3 ) \neg \
(\partial / 
\partial z_i ) (\partial \theta(z)/ \partial z_i )^{-1} $ 
($\neg $ 
is
the contraction operator).

Whence follows that the involution $ z \ra -z $ acts on the image of
$\HH^0(\hol_A(S))$  in $\HH^0(\hol_S(K_S)$ as multiplication by $-1$.

Let us now choose  in particular a surface $S$ which is the inverse
image of a theta
divisor $\Theta$ on $J$: then the subspace $V_{++}$ coming from $
\HH^0(\Omega^2_A)$
is the pull back of
$\HH^0(\Omega^2_{\Theta})$, so it consists of the sections in 
$\HH^0(\Omega^2_S) =
\HH^0(\hol_S (K_S))$ which are invariant under the fixed point free covering
involution  $z \ra z + \eta$ for  the double cover $\pi : S \ra \Theta$.

On our particular surface $S$ acts the group $(\Z/ 2)^2$ generated by
$z \ra z +
\eta$ and by $ z \ra - z $ for our choice of the origin ( c.f. Prop.
\ref{even functions}), and we
see that, if we define $V_{--}$  as the one dimensional space coming from
$\HH^0(\omega_A(S))$, then  $V_{--}$ is an eigenspace with eigenvalue
$-1$ for both
the  involutions above.

In particular, it follows that the involution $\iota$ defined by $\iota (z) =
- z + \eta$ acts trivially on the space $\HH^0(\hol_S(K_S))$.
Therefore the canonical map of such a special $S$ factors through the
involution
$\iota$.

\bigskip

\ \  \    {\bf Geometry of the situation for special surfaces }

\medskip

Let $ Z : = S/ \iota$.

\begin{lem}
The involution $ \iota$ has exactly $32$ isolated fixed points on $S$.
\end{lem}

{\it Proof of the lemma.}
Let us find the fixed points of $ \iota$ recalling that $\iota ( z)=
- z + \eta$. Then $z$ yields a fixed point on $A$ iff $ 2 z \equiv \eta \ (mod
\Lambda)$. The fixed points moreover lie on $S$ if and only if they project
in $ J = \C^3 / \Lambda(\Omega)$ to a ($2$-torsion) point which lies on
$\Theta$, i.e., to an odd thetacharacteristic.

Set $\Lambda' := \Lambda(\Omega)$, thus $\eta \in \Lambda'$, whence for such
a fixed point $2z \in \Lambda'$ and its image in $\Lambda' / \Lambda \cong \Z /
2$ is non trivial.

Therefore the number of the odd thetacharacteristics which are image of such
a fixed points are in bijection with the set $ N \subset ((\Z/ 2)^3)^2$
defined by the following equations:
$$N = \{ (x,y)| ^tx y = 1 , x_1 = 1 \}    .$$
Whence, $card(N)= 16$ and there are exactly 32 fixed points on $S$.

\qed

\begin{rmk}
Since the double cover $ S \ra Z$ is ramified exactly on the 32
corresponding nodes of $Z$, these form an even set according to the
definition of Catanese [1981].

Note moreover that $\iota$ acts as multiplication by $-1$ on the space of
global
$1-$forms, therefore the quotient surface $Z$ has $ q(Z)=0, p_g(Z)= 4, K^2_Z=
6.$
\end{rmk}

Then the canonical map of $S$, for $S$ special,
factors through $Z$. In turn, since $V_{++}$ is base point free, this
means  that there is a
point $O \in \PP^3 - \Sigma$ so that the projection with centre $O$ to $\PP^2$
yields the composition of the projection onto $\Theta$ with the
canonical map of
$\Theta$.

On the other hand, as well known, $\Theta$ is the symmetric product
$C^{(2)}$ of a
curve
$C$ of genus $3$, which, since $\Theta$ is smooth, is a smooth plane
quartic curve
$ C = C_4 \subset \PP^2$.

CLAIM: the canonical map of $C^{(2)}$ sends the divisor $ P+Q$ to the line
generated by $P$ and $Q$.

{\it Proof of the claim}.
If $\omega_1, \omega_2,\omega_3$ are a basis of $\HH^0(\Omega^1_C)$,
then a basis of
the canonical system of $C^{(2)}$ is given, on the Cartesian product $C^2$, by
$\omega_i(P) \wedge \omega_j (Q) + \omega_i(Q) \wedge \omega_j (P)$,
but this vector
is the wedge product of the two vectors $\omega_i(P)$ and $\omega_j(Q)$.

That this is a morphism follows e.g. since its base locus on $C^2$ is just the
diagonal $\Delta$, but $ \epsilon^* |K_{C^{(2)}}| =  
|K_{C^2} - \Delta |$, whence
$|K_{C^{(2)}}|$ is free from base points.

\qed

We let now $Y $ be the quotient of $ \Theta $ by the multiplication
by $-1$, whence
$ Y = S / (\Z/ 2)^2$: $Y$ has $ K_Y^2 = 3$, $ q(Y)=0, p_g(Y) = 3$ and
its canonical
map is a triple cover of $\PP^2$, branched on the dual curve $C^{\vee}$ of $C$.
In fact, multiplication by $-1$ on $\Theta$ corresponds to
residuation with respect
to $K_C$ on $C^{(2)}$.

$Y$ has $28$ nodes, corresponding to the odd thetacharacteristics of
$C$. The covering $ Z
\ra Y$ is etale, except over  $12$ of the nodes of $Y$: as we saw, $Z$ has
exactly $32$ nodes lying above the remaining $16$ nodes of $Y$, over these
$12$ nodes lie instead $12$ smooth points of $Z$.

\begin{rmk}

1) The bicanonical system of $C^{(2)}$ (cf. Catanese, Ciliberto and 
Mendes-Lopes [1998]) factors
through the
bicanonical system of $Y$, which embeds $Y$ in $\PP^6$, since it is
induced by the
sections of $ \HH^0( J, \hol_J(2 \Theta))$.

2) The monodromy of $\Theta \ra \PP^2$ is the full symmetric group
$\mathcal S_4$.
The monodromy of the canonical map of $Z$  is instead the symmetric group
$\mathcal S_3$.
\end{rmk}

\begin{lem}
The image $\Sigma$ of $Z$ is a surface of degree $6$( hence,
birational to $Z$).
\end{lem}

{\it Proof of the lemma}.
Consider the morphism $f : Z \ra \PP^2$, obtained  as the composition of the
canonical map $\phi$ of $Z$  with the projection  $p$ with centre
$O$ of $\Sigma$
to $\PP^2$.

It cannot be that $ \deg (\Sigma)=2$, otherwise $p$ would be branched 
on a plane
conic, whereas the branch curve of $f$  is the irreducible curve $C^{\vee}$.

If instead $ \deg (\Sigma)=3$, then there would be a covering 
involution $i$ for
$\phi$. Since there is already a covering involution $j$ for $f$,
gotten from the
double cover $ Z \ra Y$, we let $G$ be the group of covering involutions for
$f$. Since the canonical map of $S$ does not factor through the one
of $\Theta$,
it follows that $ i \neq j$.

Then $G$ is a group of order $h \geq 4$ with $h$ dividing $6$, thus $h=6$ and
$f$ should be Galois.

This is however a contradiction, since the inverse image of the branch curve
$C^{\vee}$ has components of multiplicity both $1$ and $2$; this holds because $Z
\ra Y$ is etale in codimension $1$, while $Y \ra \PP^2$ has simple branching on the
curve $C^{\vee}$, and the general tangent to $C^{\vee}$ is not a bitangent.

\qed

With the result of the previous lemma in our hands, we can finish the
proof of the
theorem. Assume that the canonical map of $S$ were always not birational.

Since for special $S$ the degree equals $2$, we would have that the
canonical map
always factors through the involution $\iota$. But, since $S$ always admits the
involution $ z \ra -z$, then $S$ would be stable under the involution
$ z \ra z + \eta$, i.e., would be a pull back of a theta divisor.
Contradicting that the Kuranishi family has dimension $7$ and not $6$.

Finally, in the special case, the surface $Z$ is a canonical model with
$K^2=6, p_g =4, q=0$ and with birational canonical map. Therefore, the
double curve of $\Sigma$ is a plane cubic, cf. Catanese [1984b].

\QED

In the special case, the equations of $\Sigma$ can be written explicitly. In fact,
giving an unramified double covering of a non hyperelliptic curve $C$ of genus $3$
is equivalent (cf. e.g. Catanese [1981]) to writing the equation of its canonical
model as the determinant of a $2 \times 2$ symmetric matrix of quadratic forms.

We have, more precisely, coordinates $x_0,x_1,x_2$ in $\PP^2$ and quadratic
forms
$Q_{33}(x),Q_{34}(x),Q_{44}(x)$ such that $$C = \{
(x_0,x_1,x_2)|Q_{33}(x) Q_{44}(x)- Q_{34}(x)^2 = 0 \};$$ moreover, the double
unramified covering $C'$ of $C$ is the genus $5$ curve whose canonical
model in $\PP^5$ is defined as the following intersection of three quadrics:
$$C' = \{
(x_0,x_1,x_2, y_3,y_4)| y_3^2 = Q_{33}(x),y_4^2= Q_{44}(x),y_3 y_4=
Q_{34}(x)
\}.
$$

Now, there is a natural surjection of $(C')^2$ onto $S$.  In fact, $\Theta$ is the
symmetric square of $C$, and thus dominated by $C^2$, and $S$ is the quotient of 
$(C')^2$ under the $(\Z/2)^2$ action permuting the the coordinates and acting with
the diagonal action of the involution $\iota : C' \ra C'$.

We can then read  the canonical map of
$S$  as the map corresponding to the $(\Z/2)^2$-invariant sections of $K_{(C')^2}$.

Recall that on the first curve of the product
$(C')^2$ a basis of $\HH^0 (K_{C'})^+ $ is given by $x_0,x_1,x_2$, and a basis
for $\HH^0 (K_{C'})^- $ is given by $y_3,y_4$. Similarly we have a basis $w_0,w_1,w_2 , z_3,
z_4$ for the second curve.

We find therefore that a basis for the $(\Z/2)^2$-invariant sections of
$K_{(C')^2}$ is provided by $u_0 := x_1 w_2 - x_2 w_1, u_1 := x_0 w_2 - x_2 w_0,
u_2 := x_0 w_1 -x_1 w_0  , v := y_3 z_4 - y_4 z_3 $( these are just
$\iota$-invariant  Pl\"ucker coordinates of the line spanned by the two points of
$C'$).

Let $a,b,c$ be the symmetric $ 3 \times 3$ matrices yielding the respective
quadratic forms $Q_{33}(x),Q_{34}(x),Q_{44}(x)$: then the entries of the matrix
$\alpha$ are polynomial functions in the respective entries of $a,b,c$ and in the
coordinates $(u_0, u_1, u_2, v)$ on $\PP^3$.

The shape of $\alpha^+$ is 
$$
\begin{pmatrix} 

v^5+Av^3+Bv & C \cr
C & v \cr
\end{pmatrix}$$

where for instance $A=
 ^t u(- 2\Lambda^2 b+\Lambda^2(a+c)-\Lambda^2a -\Lambda^2b)u
$,  and 
$$C= \det
\begin{pmatrix}
 ^t x a x & ^t w a x &^t w a w \\
 ^t x b x & ^t w b x &^t w b w \\
 ^t x c x & ^t w c x &^t w c w \\
\end{pmatrix}.$$

We have not yet found a compact expression for $B$, the one we have is too long to
be reproduced anywhere.

\section{ Irregular surfaces with $p_g=4, q=2$}

This section will be devoted to the description of another interesting example,
of surfaces with the following invariants: $p_g=4, q=2, K^2=18$ and
birational canonical morphism onto its image.

The surfaces are obtained as $(\Z/2\Z)^2$-Galois covers of a principally
polarized Abelian surface $A$, with branch locus consisting of
$3$ divisors $D_1, D_2, D_3$ which are algebraic equivalent to
the theta divisor $\Theta$. We shall follow the notation of Catanese [1984a].

We choose then $L_1, L_2, L_3$ divisors which are also algebraically
equivalent to $\Theta$, and such that
$$ 2 L_i \equiv D_j + D_k, \\ \forall i \neq j \neq k \neq i.$$

We take the corresponding $(\Z/2\Z)^2$-Galois cover $\pi: S \rightarrow A$
such that
$$\pi_* \hol_S = \hol_A \bigoplus (\oplus_{i=1,2,3} \hol_A
(-L_i)), \pi_* \omega_S = \hol_A \bigoplus (\oplus_{i=1,2,3} \hol_A
(L_i)) .$$
It follows immediately that the constructed surfaces have the numerical
invariants as desired: for instance, since $K_S$ is the ramification
divisor $ R$, and $2R \equiv \pi^* (D)$, where $ D = D_1+ D_2+ D_3$,
we have $K_S^2 = R^2 = D^2 = 9 \Theta^2 = 18$.

We recall the standard notation, by which $D_i = div(x_i)$, $R_i = div
(z_i)$ so that $S$ is defined by the equations 
$$w_i^2 = x_j x_k , \ 
w_i x_i = w_j w_k$$ in
the rank $3$ bundle $(\oplus_{i=1,2,3} \hol_A
(L_i)),$ and  we have $z_i^2 = 
x_i, w_i = z_j z_k$. 

We also let $\phi_i$ be the
unique section of 
$\hol_A
(L_i)) $, and $C_i := div (\phi_i)$.

With this notation, 
there are $4$ sections of the canonical sheaf
$\omega_S$, namely : $ 
\omega: = z_1 z_2 z_3$, and $\forall i=1,2,3$,
$\omega_i:= \omega 
\phi_i /w_i = z_i \phi_i. $

We obtain immediately that the base 
locus of the canonical system projects
down in $A$ to the set $ D 
\cap (\cap_{i=1,2,3} (D_i + C_i)).$

\begin{rmk}
The surface $S$ has 
base point free canonical system provided the
$6$ curves $D_i, C_i$ 
have no point common to three of them.
Since any three of the six 
curves can be chosen as arbitrary translates of
the theta divisor, it 
follows easily that for a general choice there are
no base points of 
$K_S$.
\end{rmk}

We assume henceforth the canonical system to be 
base-point free, so that we
have the canonical morphism
$$\Phi : S 
\rightarrow \Sigma \subset \PP^3$$
and we use the characters of the 
Galois group in order to study the
geometry of the map $\Phi$ and 
more generally the canonical ring of $S$.

We have here $\sR (S) = 
\bigoplus^{ \infty}_{m=0} \HH^0(\hol_S(m 
K_S))$
where
$$\HH^0(\hol_S(2m K_S)) = \HH^0(\hol_A (mD))\bigoplus 
\HH^0(\oplus_{i}
\hol_A(-L_i + m D))$$
$$\HH^0(\hol_S((2m+1) K_S)) = 
\omega \HH^0(\hol_A (mD))\bigoplus
(\oplus_{i} z_i \HH^0(
\hol_A(+L_i 
+ m D))).$$

We have $4$ generators for $\sR (S)$ in degree $1$, namely,
$\omega, \omega_1,\omega_2,\omega_3,$ moreover we observe the following
dimensions for the four respective eigenspaces in degree $2$: $ \dim (\sR
(S)_0) = 9$, $ \dim (\sR (S)_i) = 4$.

\begin{lem}
$\Phi (S)$ is not a quadric (if the canonical system is base point free).
\end{lem}
\proof
It suffices to show the linear independence of the $10$ monomials
$ \omega^2, \omega_i^2$ and $\omega \omega_i, \omega_j \omega_k$.
Any linear relation is a sum of linear relations in each eigenspace, and
clearly, if $i,j,k$ are distinct indices, then $\omega \omega_i,
\omega_j \omega_k$ are independent since their divisors are
$ 2 R_i + R_j + R_k + C'_i$, $  R_j + R_k + C'_j + C'_k$ respectively,
$C'_i$ being the inverse image of the divisor $C_i$.

Moreover, a linear relation of the form
$\Sigma_{i=0,1,2,3} c_i \omega_i^2 =0 $ would translate into a relation
$c_0 x_1 x_2 x_3 + \Sigma_{i=1,2,3} c_i x_i \phi_i^2 = 0$ and since
w.l.o.g. we may assume that $c_3=1$, we obtain that $\omega_3$ vanishes at
the points where $x_1 = x_2 = 0$, contradicting that the canonical system
is base point free.

\QED

   \begin{teo}
For general choice of the three divisors $D_i$ the canonical map  $\Phi$
is birational onto its image.
   \end{teo}
\Proof
Consider the $8$ points $P \in C'_1 \cap C'_2 = \{ \phi_1 = \phi_2 = 0\}$.
They map to $ (z_1 z_2 z_3 (P), 0,0, z_3 \phi_3 (P))$.
Moreover, the inverse image of these points is contained in $\omega_1 =
\omega_2= 0$ which consists of these $8$ points, plus points in $R_1$ or in
$R_2$ which therefore map to points where the first coordinate equals $0$.

Thus, the inverse image of the punctured line $y_1= y_2 =0, \ y_0 \neq 0$
consists of these $8$ points, which form two $(\Z/2)^2$ orbits.
For each point $ (a,0,0,b)$ in the image, since by generality we may
assume $ a \neq b \neq 0 \neq a$, the inverse image consists therefore of
either $2$ or $4$ points.
However, $ \deg (\Phi) \deg (\Sigma) = 18$, whence the only possibility
that $\Phi$ may not be birational is that $\deg (\Phi)=2$.

Assume this to be the case: then, since  $(\Z/2)^2$  acts
$\Phi$-equivariantly on $S$ and on $ \PP^3$, then we would have an
involution $i$ on $A$ which would lift to $S$, and actually in such a way
to centralize the Galois action. This however implies that
$i$ leaves the three branching divisors $D_i$ invariant.

Consider then the curve $D_1$, which has genus $2$. It possesses then only
the hyperelliptic involution, or at most a finite number of  involutions
whose quotient is an elliptic curve. Since however we may choose $D_2$
to cut $D_1$ in any assigned pair of points of $D_1$, we easily get
the desired contradiction.

\QED

\vspace{10 mm}

{\bf AUTHOR'S ADDRESS}\\

Fabrizio Catanese and

Frank-Olaf Schreyer

Lehrstuhl Mathematik VIII

Universit\" at Bayreuth

D- 95440 Bayreuth (Germany)


\begin{thebibliography}{AAAAA}



\bibitem[A-O]{1}  V. Ancona, G. Ottaviani, ``An introduction to the
                      derived categories and the theorem of
                      Beilinson", Atti Acc. Peloritana dei Pericolanti,
                         LXVII (1989), 99-110.
\bibitem[A-S]{0} E. Arbarello, E. Sernesi, ``The equation of a 
plane curve'', Duke Math. J. 46 (1979), 469-485.


\bibitem[Ash]{2}    T. Ashikaga, ``A remark on the geography of
         surfaces with birational canonical morphisms",
         Math. Ann. 290, 63-76 (1991).

\bibitem[At]{}  M.F.  Atiyah, ``Vector bundles   over an elliptic curve''
Proc. Lond. Math. Soc., III. Ser. 7, 414-452 (1957).

\bibitem[Bar]{7}    W. Barth, ``Counting singularities of quadratic forms
                 on vector bundles",in 'Vector Bundles and
                 Differential equations', Proc. Nice 1979, Birkhauser,
                 P.M. 7 (1980), 1-29.

\bibitem[B-P-V]{10}           W. Barth, C. Peters A. van de Ven,
``Compact complex surfaces'', Springer Ergebnisse F.3, vol. 4 (1984),
                 Berlin Heidelberg

\bibitem[Bau]{11}              I.C. Bauer, 
Surfaces with $K^2=7$ and $p_g=4$.
Mem. Am. Math. Soc. 721 (2001) 

\bibitem[Bea78]{12}    A. Beauville, ``Surfaces algebriques complexes",
                   Asterisque 54, Soc. Math. France, Paris (1978).

\bibitem[Bea79]{13}    A. Beauville, ``L'application canonique pour les
                 surfaces de type general", Inv. Math. 55 (1979),
                 121-140.

\bibitem[Bea82]{15}      A. Beauville, ``L'inegalit\'e $p_g\geq 2q-4$
                          pour les
                         surfaces de type general", Bull.Soc.Math. France
                         110 (1982), 344-346.

\bibitem[Bei]{16}       A. Beilinson, ``Coherent sheaves on $\PP^n$ and
                         problems of linear algebra", Funct.Anal. Appl.
                         12 (1978), 214-216 ( translated from
                         Funkt.Anal. i Priloz. 12, 3 (1978) 68-69).
\bibitem[B-G-G]{17}
I.N. Bernshtein, I.M. Gel'fand, S.I. Gel'fand,"Algebraic
                         bundles on $\Pi^n$ and problems of linear algebra",
                         Funct.Anal. Appl. 12 (1978), 212-214
                         ( translated from Funkt.Anal. i Priloz. 12, 3
                         (1978) 66-68).

\bibitem[Bol]{}              O. Bolza, ``On binary sextics with linear transformations into
themselves", Amer. Jour. Math. 10 (1888), 47-
         70.

\bibitem[Bom]{18}              E. Bombieri, ``Canonical models of surfaces
of         general type", I.H.E.S. Publ.Math. 42 (1973), 171-
         219.



\bibitem[B-E]{23}    D.A. Buchsbaum, D. Eisenbud, ``What annihilates
         a module?"       J. Algebra 47 (1977), 231-243.


\bibitem[Can00]{} A. Canonaco,
``A Beilinson-type theorem for coherent sheaves on weighted projective spaces'',
J. Algebra 225 (2000), 28-46.

\bibitem[Can02]{25}     A. Canonaco, ``
         ", Thesis, Scuola Normale Superiore Pisa (2002)


\bibitem[C-E]{28}     G. Casnati, T.Ekedahl, ``Covers of algebraic
         varieties I. General structure theorem, covers of
         degree 3,4 and Enriques surfaces", J.Alg.Geom.5 (1996), 439-460.

\bibitem[Cas91]{29}     G. Castelnuovo, ``Osservazioni intorno alla
         geometria sopra una superficie, I, II ``, Rendiconti
         del R. Istituto Lombardo, s. II, 24 (1891), also in
         'Memorie scelte ', Zanichelli  (1937), Bologna,
         245-265.

\bibitem[Cas05]{30}    G. Castelnuovo, ``Sulle superficie aventi il genere
         aritmetico negativo", Rend. Circ. Mat. Palermo, 20
         (1905), 55-60.

\bibitem[Cat81]{31}    F. Catanese, ``Babbage's conjecture, contact of
                 surfaces, symmetric determinantal varieties
                 and applications", Inv. Math. 63 (1981), 433-465.

\bibitem[Cat84a]{32}    F. Catanese, ``On the moduli space of surfaces of
         general type", J. Differential Geom., 19 (1984),
         483-515.

\bibitem[Cat84b]{33}    F. Catanese, ``Commutative algebra methods and
         equations of regular surfaces", in 'Algebraic
         Geometry - Bucharest' 1982, Springer L.N.M.,
         1056 (1984), 68-111.

\bibitem[Cat85]{34}    F. Catanese, ``Equations of pluriregular varieties
         of general type", in ``Geometry today-Roma
         1984", Progr. in Math. 60, Birkhauser (1985), 47-
         67.


\bibitem[Cat97]{38}    F. Catanese,
``Homological algebra and algebraic surfaces''
J. Kollar (ed.) et al., ``Algebraic geometry''
Proceedings of the Summer Research Institute, Santa Cruz, CA, USA, 
July 9-29, 1995.
Providence, RI: American Mathematical Society. Proc. Symp. Pure Math. 
62 (1997), 3-56.


\bibitem [C-F]{} F. Catanese, M. Franciosi,
``Divisors of small genus on algebraic surfaces and projective embeddings,''
Teicher, Mina (ed.), Proceedings of the Hirzebruch 65 conference on 
algebraic geometry,
  Bar-Ilan University, Ramat Gan, Israel, May 2-7, 1993. Ramat-Gan: 
Bar-Ilan University,
Isr. Math. Conf. Proc. 9 (1996), 109-140.

\bibitem[C-F-H-R]{} F. Catanese, M. Franciosi, K. Hulek, 
M. Reid,
``Embeddings of curves and surfaces'',
Nagoya Math. J. 154 (1999), 185-220.


\bibitem[C-C91]{} F. Catanese, C. Ciliberto,
``Surfaces with $p\sb g = q = 1$'', 
F. Catanese et al. (eds), ``Problems in the theory of surfaces and 
their classification,''
  Papers from the meeting held at the Scuola Normale Superiore, 
Cortona, Italy, October
10-15, 1988. London: Academic Press. Symp. Math. 32 (1991), 49-79.


\bibitem[C-C93]{} F. Catanese, C. Ciliberto,
``Symmetric products of elliptic curves and surfaces of general type 
with $p\sb g=q=1$'',
J. Algebr. Geom. 2, (1993), 389-411.

\bibitem[C-C-ML]{} F. Catanese, C. Ciliberto, Mendes-Lopes,
`` On the classification of irregular surfaces of general type with
  nonbirational bicanonical map,''
Trans. Am. Math. Soc. 350, No.1, 275-308 (1998).


\bibitem[Cil81]{44}            C. Ciliberto, ``Canonical surfaces with
$p_g=p_a=4$ and
     $K^2 = 5,\dots,10$", Duke Math. J. 48 (1981), 121-157.

\bibitem[Cil83]{46}            C. Ciliberto, ``Sul grado dei generatori dell'
anello                canonico di una superficie di tipo generale", Rend.
                 Sem. Mat. Torino 41, 3, (1983), 83-111.

\bibitem[Clem00]{46}  H. Clemens, ``Cohomology and obstructions I: on the geometry of formal
Kuranishi theory``, math.AG/9901084 v2.


\bibitem[De]{48}              O. Debarre, ``Inegalites numeriques pour les
surfaces
                 de type general", Bull. Soc. Math. France 110 (1982),
                 319-346.

\bibitem[Di]{50}              A. C. Dixon, ``Note on the reduction of a
ternary                quantic to a symmetrical determinant", Proc.
                 Camb. Phil. Soc. 11 (1902), 350-351.

\bibitem[Ei80]{}  D. Eisenbud, ``Homological algebra on
   a complete intersection'', Trans. Amer. Math. Soc. 260 (980), 35-64.


\bibitem[Ei95]{51}              D. Eisenbud, ``Commutative Algebra with a
                   view towards Algebraic Geometry ``, Springer
                 G.T.M. 150, New York (1995) .

\bibitem[E-F-S]{} D. Eisenbud, G. Fl\o ystad, F.-O. Schreyer,
   ``Sheaf cohomology and free fesolutions over exterior algebras'',
    preprint,  math.AG/0104203.

\bibitem[E-S]{} D. Eisenbud, F.-O. Schreyer,
``Resultants and Chow forms via exterior syzygies'', preprint,  math.AG/0111040


\bibitem[En]{54}              F.Enriques, `` Le Superficie Algebriche",
Zanichelli,                        Bologna (1949)

\bibitem[Fit]{56}             H. Fitting, ``Die Determinantenideale einer
                           Moduls", Jahresber. Deutscher Math. Verein. 46
                         (1936), 195-220.

\bibitem[Fuj]{} T. Fujita, 
``On K\"ahler fiber spaces over curves'',
J. Math. Soc. Japan 30, 779-794 (1978)

\bibitem[Gra]{58}             M. Grassi, `` Koszul modules and Gorenstein
                           algebras", J. Alg. 180, (1996) 918-953.

\bibitem[Har]{66}    R. Hartshorne, Algebraic Geometry, Springer G.T.M.
                 52, New York (1977)

\bibitem[Hilb]{67}            D. Hilbert, ``\"Uber die Theorie der
Algebraischen Formen", Math. Ann. 36 (1890), 473-534.

\bibitem[Hor75]{69}            E. Horikawa, ``On deformation of quintic
surfaces",                 Inv. Math. 31 (1975),43-85.
\bibitem[Hor1-5]{70}          E. Horikawa, ``Algebraic surfaces of general
type withsmall $c_1^2$", I, Ann. of Math. (2) 104 (1976), 357-387;    
II, Inv. Math. 37 (1976), 121-155; III, Inv. Math. 47 (1978), 209-248; 
IV, Inv. Math. 50 (1979), 103-128; V, J. Fac. Sci. Univ. Tokyo, Sect. A. Math.
283 (1981), 745-755.

\bibitem[dJ-vS]{72}   T. de Jong, D. van Straten, ``Deformations of the
                 normalization of hypersurfaces", Math. Ann. 288
                 (1990), 527-547
\bibitem[J-L-P]{73}   T. Josefiak, A. Lascoux, P. Pragacz, ``Classes of
                 determinantal varieties associated with symmetric
                 and skew-symmetric matrices", Izv. An. SSSR 45,9
                 (1981), 662-673, translated in Math.USSR Izv. 18
                 (1982), 575-586.

\bibitem[Kap]{74}     M. Kapranov, ``On the derived categories of coherent
                 sheaves on some homogeneous spaces", Inv. Math. 92
                 (1988), 479-508.

\bibitem[Kod]{77}             K. Kodaira, `` On characteristic systems of
families of surfaces with ordinary singularities in a projective space", 
Amer. J. Math. 87 (1965), 227-256 .

\bibitem[Kon91]{79}             K. Konno, ``A note on surfaces with pencils of
non-hyperelliptic curves of genus 3'', Osaka J. Math. 28 (1991), 737-745.

\bibitem[Kon96]{80}             K. Konno, ``A lower bound of the slope of
trigonal fibrations", Int. J. Math. 7, 1, (1996) 19-27.

\bibitem [Kon93]{} K. Konno, 
`` Non-hyperelliptic fibrations of small genus and certain irregular canonical surfaces", 
Ann. Sc. Norm. Super. Pisa, Cl. Sci., IV. Ser. 20, No.4, 575-595 (1993)

\bibitem[Miy]{87}  Y. Miyaoka, ``On the Chern numbers of surfaces
of general type", Inv. Math. 42 (1977), 225-237.

\bibitem[M-P]{88}   D. Mond, R. Pellikaan,
``Fitting ideals and multiple  points of analytic mappings",
Springer L.N.M. 1414,New York-Berlin (1987), 54 pages

\bibitem[Mum]{} D. Mumford, 
`` On the equations defining Abelian varieties. I-III ",
Invent. Math. 1, 287-354 (1966); ibid. 3, 75-135, 215-244 (1967).


\bibitem[Ran] { } Z. Ran, 
`` Hodge theory and deformations of maps. ``
Compos. Math. 97, No.3, 309-328 (1995).

\bibitem[Rao]{100}            P. Rao, ``Liaison among curves in $\PP^3$",
Inv.Math.50                 (1979), 205-217.

\bibitem[Reid]{103}           M. Reid, Math. Review 86c: 14027.
\bibitem[Reid79]{104}   M. Reid, ``$\pi_1$ for surfaces with small  $c_1^2"$,
in         'Algebraic Geometry', Springer LNM 732, 534-544
         (1979).

\bibitem[Sch]{107}    F.-O. Schreyer, `` Small fields and constructive
         algebraic geometry",n: M. Maruyama, editor: Moduli of vector
      bundles, Marcel Dekker Inc., New York 1996, 221-228

\bibitem[Sern]{108}    E. Sernesi, ``L'unirazionalita' della varieta' dei
moduli         delle curve di genere dodici", Ann. Scuola Norm. Pisa
         8 (1981), 405-439.

\bibitem[Ser]{109}    J.P. Serre, ``Faisceaux algebriques coherents", Annals
         of Math. 61, 2 (1955), 197-278 .

\bibitem[Sim]{} C. Simpson, 
``Subspaces of moduli spaces of rank one local systems.''
Ann. Sci. Éc. Norm. Supér., IV. Sér. 26, No.3, 361-401 (1993)



\bibitem[Taka96]{} T. Takahashi,
``Certain algebraic surfaces of general type with irregularity one and
their canonical mappings.''
Tohoku Mathematical Publications. 2. Sendai: Tohoku Univ., 
Mathematical Institute, vi,
60 p. (1996).

\bibitem[Taka98]{} T. Takahashi, ``
Certain algebraic surfaces of general type with irregularity one
and their canonical mappings.''
Tohoku Math. J., II. Ser. 50, No.2, 261-290 (1998).

\bibitem[Xiao87]{}  Xiao,Gang
`` Fibered algebraic surfaces with low slope.'' 
Math. Ann. 276, 449-466 (1987)

\bibitem[Yau]{119}     S. T. Yau, ``Calabi's conjecture and some new
         results in algebraic geometry", Proc. Nat. Acad. Sc.
         USA 74(1977), 1798-1799.

\bibitem[Zuc94]{121}      F. Zucconi, `` Su alcune questioni relative alle superficie di tipo
generale con mappa canonica composta con un fascio o di grado $3$``, Tesi di Dottorato,
Universita' di Pisa, 1994
\end{thebibliography}
\end{document}